\documentclass[11pt]{article}
\usepackage{amssymb}
\usepackage{amsmath}
\usepackage[all]{xy}
\usepackage{times}
\usepackage{graphicx}

\title{Some complete $\omega$-powers of a one-counter language,\\ for any Borel class of finite rank}
\author{Olivier FINKEL and Dominique LECOMTE$^1$}
\date{\today}

\def\ufootnote#1{\let\savedthfn\thefootnote\let\thefootnote\relax
\footnote{#1}\let\thefootnote\savedthfn\addtocounter{footnote}{-1}}

\newcommand{\fa}{\forall}
\newcommand{\Ga}{\Gamma}
\newcommand{\Gas}{\Gamma^{<\om}}
\newcommand{\Si}{\Sigma}
\newcommand{\Sis}{\Sigma^{<\om}}
\newcommand{\ra}{\rightarrow}

\newcommand{\la}{language}
\newcommand{\ite}{\item}

\newcommand{\ol}{ $\omega$-language}
\newcommand{\orl}{ $\omega$-regular language}
\newcommand{\om}{\omega}

\newcommand{\noi}{\noindent}
\newcommand{\tla}{\twoheadleftarrow}
\newcommand{\de}{deterministic }
\newcommand{\proo}{\noi {\bf Proof.} }
\newcommand {\ep}{\hfill $\square$}
\newcommand{\Sio}{\Si^\om}

\newcommand{\borom}{{\bf\Delta}^{0}_{\omega}}
\newcommand{\borxi}{{\bf\Delta}^{0}_{\xi}}
\newcommand{\bormxi}{{\bf\Pi}^{0}_{\xi}}

\newcommand{\bormone}{{\bf\Pi}^{0}_{1}}

\newcommand{\bormtwo}{{\bf\Pi}^{0}_{2}}

\newcommand{\bormom}{{\bf\Pi}^{0}_{\omega}}

\newcommand{\borel}{{\bf\Delta}^{1}_{1}}

\newcommand{\borone}{{\bf\Delta}^{0}_{1}}
\newcommand{\bortwo}{{\bf\Delta}^{0}_{2}}
\newcommand{\borthree}{{\bf\Delta}^{0}_{3}}
\newcommand{\boraone}{{\bf\Sigma}^{0}_{1}}
\newcommand{\boratwo}{{\bf\Sigma}^{0}_{2}}
\newcommand{\borathree}{{\bf\Sigma}^{0}_{3}}
\newcommand{\boraom}{{\bf\Sigma}^{0}_{\omega}}
\newcommand{\boraxi}{{\bf\Sigma}^{0}_{\xi}}
\newcommand{\ana}{{\bf\Sigma}^{1}_{1}}

\newtheorem{thm} {Theorem}
\newtheorem{defi} [thm] {Definition}
\newtheorem{cor} [thm] {Corollary}
\newtheorem{lem} [thm] {Lemma}
\newtheorem{prop} [thm] {Proposition}

\begin{document}

\maketitle

\centerline{$\bullet$ CNRS, Universit\'e Paris Diderot, Sorbonne Universit\'e,} 

\centerline{Institut de Math\'ematiques de Jussieu-Paris Rive Gauche, Equipe de Logique Math\'ematique}

\centerline{Campus des Grands Moulins, b\^atiment Sophie-Germain, case 7012, 75205 Paris cedex 13, France}

\centerline{finkel@math.univ-paris-diderot.fr}\medskip

\centerline{$\bullet^1$ Sorbonne Universit\'e, Universit\'e Paris Diderot, CNRS,} 

\centerline{Institut de Math\'ematiques de Jussieu-Paris Rive Gauche, Equipe d'Analyse Fonctionnelle}

\centerline{Campus Pierre et Marie Curie, case 247, 4, place Jussieu, 75 252 Paris cedex 5, France}

\centerline{dominique.lecomte@upmc.fr}\medskip

\centerline{$\bullet^1$ Universit\'e de Picardie, I.U.T. de l'Oise, site de Creil,}

\centerline{13, all\'ee de la fa\"\i encerie, 60 107 Creil, France}\medskip\medskip\medskip\medskip\medskip\medskip

\ufootnote{{\it 2010 Mathematics Subject Classification.}~Primary: 03E15, Secondary: 54H05, 68R15}

\ufootnote{{\it Keywords and phrases.}~Borel class, complete, context-free, one-counter automaton; $\omega$-power}

\ufootnote{{\it Acknowledgements.}~We thank very much the anonymous referees for their very useful comments about a preliminary version of our article.}

\noindent {\bf Abstract.} We prove that, for any natural number $n\!\geq\! 1$, we can find a finite alphabet $\Si$ and a finitary language $L$ over $\Si$ accepted by a one-counter automaton, such that the $\omega$-power 
$$L^\infty\! :=\!\{ w_0w_1\ldots\!\in\!\Si^\om\mid\fa i\!\in\!\omega  ~~w_i\!\in\! L\}$$ 
is ${\bf\Pi}^0_n$-complete. We prove a similar result for the class ${\bf\Sigma}^0_n$.

\newpage

\section{$\!\!\!\!\!\!$ Introduction}\indent

 We pursue in this paper the study of  the topological complexity of $\om$-powers of languages of finite words over a finite alphabet 
$\Sigma$. A finitary language over a finite alphabet $\Sigma$ is a subset $A$ of the set $\Sigma^{<\omega}$ of finite words with letters in $\Sigma$. The set of infinite words over the alphabet $\Si$, i.e., of sequences of length $\om$ of letters of   $\Si$, is denoted 
$\Si^\om$. The $\omega$-{\bf power} associated with $A\subseteq \Sigma^{<\omega}$ is the set $A^\infty$ of the infinite words constructible with $A$ by concatenation, i.e., 
${A^\infty\! :=\!\{\ a_{0}a_{1}\ldots\!\in\!\Sigma^\omega\mid\forall i\!\in\!\omega~\ a_{i}\!\in\! A\ \}}$. Notice that we denote here $A^\infty$ the $\omega$-power associated with $A$, as in \cite{Lecomte05,Fink-Lec2}, while it is often denoted $A^\om$ in Theoretical Computer Science papers, as in \cite{Staiger97,Fin01a,Fin03a,Fin-Lec}. Here we reserved the notation $A^\om$ to denote the Cartesian product of countably many copies of $A$ since this will be often used in this paper.\medskip
      
     In the theory of formal languages of infinite words, accepted by various kinds of automata, 
the $\omega$-powers appear very naturally  in the characterization of the class 
$REG_\om$  of  \orl s (respectively,  of the class $CF_\om$ of context free \ol s)  as the $\om$-Kleene closure 
of the family $REG$ of regular finitary languages (respectively,   of the family $CF$ of context free finitary languages) \cite{Staiger97}.
Since the set $\Sio$ of infinite words over a finite alphabet $\Si$ can be   equipped 
with the usual Cantor topology, the question of  the topological  complexity of  $\om$-powers of 
finitary  languages naturally arose and was  posed in particular by 
Niwinski \cite{Niwinski90},  Simonnet \cite{Simonnet92},  and  Staiger \cite{Staiger97}.  Moreover the $\om$-powers have also been studied from the perspective of Descriptive Set Theory 
in \cite{Lecomte05,Fink-Lec2}.\medskip

As the concatenation map, from $A^\omega$ onto $A^\infty$, which associates to a given sequence $(a_{i})_{i\in\omega}$ of finite words the concatened word $a_{0}a_{1}\ldots$, is continuous, an $\omega$-power is always an analytic set.\medskip

        It  was   proved in \cite{Fin03a}    that there exists a 
(context-free) language $L$ such that $L^\infty$ is analytic but not Borel. Amazingly, 
the language $L$ is  very simple to describe and it is accepted by a simple $1$-counter automaton. Louveau has proved independently 
that analytic-complete $\omega$-powers exist, but the existence was proved in a non effective way (this is non-published work).  
We refer the reader to \cite{ABB96} for basic notions about context-free languages.\medskip

Concerning Borel $\om$-powers, it was proved that,
for each integer $n\geq 1$, there exist some $\om$-powers of  (context-free) languages 
which are ${\bf \Pi}_n^0$-complete Borel sets,  \cite{Fin01a}.   It was proved in \cite{Fin04-FI} that  there exists a finitary language $V$ 
such that $V^\infty$ is a Borel set of infinite rank, and in \cite{Fin-Dup06}
that there is a (context-free) language $W$ such that $W^\infty$ is Borel above 
${\bf \Delta_\omega^0}$.\medskip

We  proved  in \cite{Fin-Lec, Fink-Lec2}  a result which showed that $\omega$-powers exhibit a great topological complexity: 
 for each nonzero countable ordinal  $\xi$, there exist $\bormxi$-complete $\omega$-powers, and ${\bf \Sigma_\xi^0}$-complete  $\omega$-powers.  This result has an effective aspect: for each recursive ordinal $\xi < \omega_1^{\rm CK}$, where  $\omega_1^{\rm CK}$ is the first non-recursive ordinal, there exists recursive finitary languages $P$ and $S$ such that $P^\infty$ is 
 $\bormxi$-complete and $S^\infty$ is ${\bf \Sigma_\xi^0}$-complete.\medskip

Many questions are still open about the topological complexity of $\om$-powers of languages in a given class like the class of context-free languages, $1$-counter languages, recursive languages, or more generally  languages accepted by some kind of automata over finite words.\medskip

In this paper we obtain the following new results about $\om$-powers of languages accepted by $1$-counter automata. 

\newpage 

\begin{thm} \label{main} Let $n\!\geq\! 1$ be a natural number.\smallskip

(a) There is a finitary language $P_n$ which is accepted by a one-counter automaton and such that the $\om$-power $P_n^\infty$ is ${\bf\Pi}_n^0$-complete.\smallskip

(b) There is a finitary language $S_n$ which is accepted by a one-counter automaton and such that the $\om$-power $S_n^\infty$ is ${\bf\Sigma}_n^0$-complete.\smallskip

Moreover, for any given integer $n\!\geq\! 1$, one can effectively construct some one-counter automata accepting such finitary languages $P_n$ and $S_n$.\end{thm}

 This article is organized as follows. Notions of automata and formal language theory are recalled in Section 2. Some basic notions of topology are recalled in Section 3. The definition and some properties of the operation of exponentiation of sets are given in Section 4. Our results related to the classes ${\bf\Pi}_n^0$ are proved in Section 5 and our results related to the classes ${\bf\Si}_n^0$ are proved in Section 6.\medskip

 We give in this article a construction of complete $\omega$-powers of a one-counter language, for any Borel class of finite rank. It remains open to determine completely the topological complexity of  $\omega$-powers of one-counter languages. Recall that it has been proved in \cite{Fin-mscs06} that for each recursive ordinal $\xi\! <\!\omega_1^{\rm CK}$, there exist some $\om$-languages $P_\xi$ and $S_\xi$  accepted by B\"uchi one-counter automata such that $P_\xi$ is $\bormxi$-complete and $S_\xi$ is ${\bf \Sigma_\xi^0}$-complete.\medskip
  
 Moreover each $\om$-language $L\subseteq \Sio$ accepted by a   B\"uchi one-counter automaton is of the form $L = \bigcup_{1\leq j \leq n} U_j\cdot V_j^\infty$, for some one-counter finitary languages $U_j$ and $V_j$, $1\leq j \leq n$. Therefore it seems plausible that there exist complete $\omega$-powers of a one-counter language, for each Borel class of recursive rank. 

\section{$\!\!\!\!\!\!$ Automata}\label{section-automata}\indent

 We assume the reader to be familiar with formal \la s, see for example \cite{HopcroftMotwaniUllman2001,Thomas90}. We first recall some of the definitions and results concerning pushdown automata and context free  \la s, as presented in \cite{ABB96,cg,Staiger97}.\medskip

 When $\Si$ is a finite alphabet, a nonempty {\bf finite word} over $\Si$ is a sequence $w\! =\! a_0\ldots a_{l-1}$, where 
$a_i\!\in\!\Sigma$ for each $i\! <\! l$, and $l\geq 1$ is a natural number. The {\bf length} of $w$ is $l$, denoted by $|w|$. If 
$|w|\! =\! 0$, then $w$ is the {\bf empty word}, denoted by $\lambda$. When $w$ is a finite word over $\Sigma$, we write 
$w\! =\! w(0)w(1)\ldots w(l\! -\! 1)$, and the prefix $w(0)w(1)\ldots w(i\! -\! 1)$ of $w$ of length $i$ is denoted by 
$w\vert i$, for any $i\!\leq\! l$. We also write $u\!\subseteq\! v$ when the word $u$ is a prefix of the finite word $v$. The set of finite words over $\Si$ is denoted by $\Si^{<\om}$, and $\Si^+$ is the set of nonempty finite words over $\Sigma$. A {\bf language} over $\Sigma$ is a subset of $\Sis$. For $L\!\subseteq\!\Sis$, the {\bf complement} 
$\Sis\!\setminus\! L$ of $L$ (in $\Sis$) is denoted by $L^-$.\medskip

  The first infinite ordinal is $\om$. An $\om$-{\bf word} over $\Si$ is an $\om$-sequence $a_0a_1\ldots$, where 
$a_i\!\in\!\Sigma$ for each natural number $i$. When $\sigma$ is an $\om$-word over $\Si$, we write 
$\sigma =\sigma(0)\sigma(1)\ldots$, and the prefix $\sigma(0)\sigma(1)\ldots\sigma(i\! -\! 1)$ of $\sigma$ of length $i$ is denoted by $\sigma\vert i$, for any natural number $i$. We also write $u\!\subseteq\!\sigma$ when the finite word $u$ is a prefix of the $\om$-word $\sigma$. The set of $\om$-words over $\Si$ is denoted by $\Si^\om$. An 
$\om$-{\bf language} over $\Sigma$ is a subset of  $\Si^\om$. For $A\subseteq \Si^\om$, the complement 
$\Si^\om\!\setminus\! A$ of $A$ is denoted by $A^-$.\medskip

 The usual {\bf concatenation} product of two finite words $u$ and $v$ is denoted $u^\frown v$ (and sometimes just $uv$). This product is extended to the product of a finite word $u$ and an $\om$-word $\sigma$: the infinite word $u^\frown\sigma$ is then the $\om$-word such that $(u^\frown\sigma )(k)\! =\! u(k)$ if $k\! <\! |u|$, and $(u^\frown\sigma )(k)\! =\!\sigma (k\! -\! |u|)$ if $k\!\geq\! |u|$.

\vfill\eject

 If $E$ is a set, $l\!\in\!\omega$ and $(e_i)_{i<l}\!\in\! E^l$, then ${^\frown}_{i<l}\ e_i$ is the concatenation 
$e_0\ldots e_{l-1}$. Similarly, ${^\frown}_{i\in\omega}\ e_i$ is the concatenation $e_0e_1\ldots$ For $L\subseteq\Sis$, 
$L^\infty\! :=\!\{ \sigma =w_0w_1\ldots\!\in\!\Si^\om\mid\fa i\!\in\!\omega  ~~w_i\!\in\! L\}$ is the $\om$-{\bf power} of $L$.

\begin{defi} 
A {\bf pushdown automaton} is a 7-tuple $\mathcal{A}\! =\! (Q,\Si ,\Ga ,q_0,Z_0,\delta ,F)$, where $Q$ is a finite set of states, 
$\Sigma$ is a finite input alphabet, $\Gamma$ is a finite pushdown alphabet, $q_0\!\in\! Q$ is the initial state, $Z_0\!\in\!\Ga$ is the start symbol which is the bottom symbol and always remains at the bottom of the pushdown stack, $\delta$ is a map from $Q\!\times\! (\Si\cup\{\lambda\} )\!\times\!\Ga$ into the set of finite subsets of $Q\!\times\!\Gas$, and $F\!\subseteq\! Q$ is the set of final states. The automaton $\mathcal{A}$ is said to be {\bf real-time} if there is no $\lambda$-transition, i.e., if $\delta$ is a map from $Q\!\times\!\Si\!\times\!\Ga$ into the set of finite subsets of $Q\!\times\!\Gas$.\smallskip

 If $\gamma\in\Ga^+$ describes the pushdown stack content, then the leftmost symbol will be assumed to be on the ``top" of the stack. A {\bf configuration} of the pushdown automaton $\mathcal{A}$ is a pair $(q,\gamma )$, where $q\!\in\! Q$ and 
$\gamma\in\Gas$. For $a\!\in\!\Si\cup\{\lambda\}$, $\gamma ,\beta\!\in\!\Gas$ and $Z\!\in\!\Ga$, if $(p,\beta )$ is in $\delta (q,a,Z)$, then we write $a\! :\! (q,Z\gamma )\!\mapsto_{\mathcal{A}}\! (p,\beta\gamma )$.\smallskip
  
 Let $w\! =\! a_0\ldots a_{l-1}$ be a  finite word over $\Si$. A  sequence of configurations $r\! =\! (q_i,\gamma_i)_{i<N}$ is called a  {\bf run of $\mathcal{A}$ on $w$ starting in the configuration} $(p,\gamma )$ if
\begin{enumerate}
\ite[(1)] $(q_0,\gamma_0)\! =\! (p,\gamma)$,

\ite[(2)]  for each $i\! <\! N\! -\! 1$, there exists $b_i\in\Si\cup\{\lambda\}$ satisfying 
$b_i\! :\! (q_i,\gamma_i)\!\mapsto_{\mathcal{A}}\! (q_{i+1},\gamma_{i+1})$ such that $a_0\ldots a_{l-1}\! =\! b_0\ldots b_{N-2}$.
\end{enumerate}

 A run $r$ of $\mathcal{A}$ on $w$ starting in configuration $(q_0,Z_0)$ will be simply called a {\bf run of} $\mathcal{A}$ {\bf on} 
$w$. The run is {\bf accepting} if it ends in a final state.\smallskip

 The language $L(\mathcal{A})$ accepted by $\mathcal{A}$ is the set of words admitting an accepting run by $\mathcal{A}$. A 
{\bf context-free language} is a finitary language which is accepted by a pushdown automaton. We denote by $CFL$ the class of context-free languages.\smallskip

 A {\bf one-counter automaton} is a pushdown automaton with a pushdown alphabet of the form 
${\Ga\! =\! \{ Z_0,z\}}$, where $Z_0$ is the bottom symbol and always remains at the bottom of the pushdown stack. A 
{\bf one-counter language} is a (finitary) language which is accepted by a one-counter automaton. 
\end{defi}

\noindent\bf Remarks.\rm\ (1) The pushdown automaton defined above is in general {\it non-deterministic}.  In the sequel, we often indicate when the considered automata can be deterministic or when the non-determinism is essential  in the behaviour of the automata.\medskip

\noindent (2) The accepting condition here is by final states. Some other accepting conditions have been considered. For instance a language is context-free if and only if it is accepted by a pushdown automaton by final states {\it and} empty stack \cite{ABB96}.  In particular, in the last sections of the paper, we will consider acceptance by final states {\it and} empty stack.

\begin{defi} Let $\Sigma ,\Gamma$ be finite alphabets.\smallskip

(a) A $(\Sigma ,\Gamma )$-{\bf substitution} is a map $f\! :\!\Sigma\!\rightarrow\! 2^{\Gamma^{<\omega}}$.\smallskip

(b) We extend this map to $\Sigma^{<\omega}$ be setting 
$f({^\frown}_{i<l}\ a_i)\! :=\!\{ {^\frown}_{i<l}\ w_i\mid\forall i\! <\! l~~w_i\!\in\! f(a_i)\}$, where 
$l\!\in\!\omega$ and $a_0,\cdots ,a_{l-1}\!\in\!\Sigma$.\smallskip

(c) We further extend this map to $2^{\Sigma^{<\omega}}$ by setting $f(L)\! :=\!\bigcup_{w\in L}~f(w)$.\smallskip

(d) Let $f$ be a $(\Sigma ,\Gamma )$-substitution, and $\mathcal{F}$ be a family of languages. If the language $f(a)$ belongs to 
$\mathcal{F}$ for each $a\!\in\!\Si$, then the substitution $f$ is called a $\mathcal{F}$-{\bf substitution}.\smallskip

(e) We then define the operation $\square$ on families of languages. Let $\mathcal{E}$, $\mathcal{F}$ be families of (finitary) languages. Then 
$\mathcal{E}~\square ~\mathcal{F}\! :=\!\{ f(L)\mid L\!\in\!\mathcal{E}\mbox{ and }f\mbox{ is a }\mathcal{F}\mbox{-substitution}\}$. 
\end{defi}

  The operation of substitution gives rise to an infinite hierarchy of context free finitary languages defined as follows.

\begin{defi} Let  $OCL(0)=REG$ be the class of regular languages,   $OCL(1)=OCL$  be  the class of one-counter languages, 
and $OCL(k\! +\! 1)\! =\! OCL(k)~\square ~OCL$, for $k\!\geq\! 1$.\end{defi}

 It is well known that the hierarchy given by the families of languages $OCL(k)$ is strictly increasing. And there is a characterization  of these languages by means of automata.

\begin{prop}[\cite{ABB96}]  A language $L$ is in $OCL(k)$ if and only if $L$ is recognized by a pushdown automaton such that, during any computation, the words in the pushdown stack remain in a bounded language of the form 
$(z_{k-1})^{<\om}\ldots (z_0)^{<\om}Z_0$, where $\{ Z_0, z_0,\ldots ,z_{k-1}\}$ is the pushdown alphabet. Such an automaton is called a $k$-{\bf iterated counter automaton}. The union $ICL\! :=\!\bigcup_{k\geq 1}~OCL(k)$ is called the family of {\bf iterated counter languages}, which is the closure under substitution of the family $OCL$.\end{prop}

 Note that we can consider that  a  $k$-iterated counter automaton is a $k$-counter automaton in the following way.  If the content of the pushdown stack of a $k$-iterated counter automaton is equal to $(z_{k-1})^{n_{k-1}} \ldots (z_0)^{n_0}Z_0$ for some natural numbers $n_0,\ldots ,n_{k-1}$, then the numbers $n_0,\ldots ,n_{k-1}$ are the contents of the counters 1,\ldots , $k$ of the $k$-counter automaton. Moreover, it is then clear that the content of the $i^{\mbox{th}}$ counter can only be changed when the contents of counters numbered $i\! +\! 1$, $i\! +\! 2$, \ldots, $k\! -\! 1$ are equal to zero. We now recall the formal definition of a $k$-counter automaton.

\begin{defi} Let $k\!\geq\! 1$ be an integer. A $k$-{\bf counter automaton} is a 5-tuple 
$\mathcal{A}\! =\! (Q,\Si ,q_0,\Delta ,F)$,  where $Q$ is a finite set of states, $\Sigma$ is a finite input alphabet, $q_0\!\in\! Q$ is the initial state,  
$$\Delta\!\subseteq\! Q\!\times\! (\Si\cup\{\lambda\} )\!\times\!\{0, 1\}^k\!\times\! Q\!\times\!\{0,1,-1\}^k$$ 
is the transition relation, and $F\!\subseteq\! Q$ is the set of final states. The $k$-counter automaton $\mathcal{A}$ is said to be 
{\bf real-time} if there is no $\lambda$-transition, i.e., if 
$\Delta\!\subseteq\! Q\!\times\!\Si\!\times\!\{ 0,1\}^k\!\times\! Q\!\times\!\{0,1,-1\}^k$.\smallskip

 If the machine $\mathcal{A}$ is in the state $q$ and $c_i\!\in\!\omega$ is the content of the $i^{\mbox{th}}$ counter 
$\mathcal{C}_i$, then the {\bf configuration} (or global state) of $\mathcal{A}$ is the $(k\! +\! 1)$-tuple $(q,c_0,\ldots ,c_{k-1})$.\smallskip

 Let $a\in \Si \cup \{\lambda\}$, $q, q'\!\in\! Q$, $(c_0,\ldots ,c_{k-1})\!\in\!\omega^k$. We write 
$$a\! :\! (q,c_0,\ldots ,c_{k-1})\!\mapsto_{\mathcal{A}}\! (q',c_0\! +\! l_0,\ldots ,c_{k-1}\! +\! l_{k-1})$$
when $(q,a,i_0,\ldots ,i_{k-1},q', l_0,\ldots ,l_{k-1})\!\in\!\Delta$, where $i_j\! =\! 0$ if $c_j\! =\! 0$ and $i_j\! =\! 1$ if $c_j> 0$. This implies that the transition relation has the property that if $(q,a,i_0,\ldots ,i_{k-1},q',l_0,\ldots ,l_{k-1})\!\in\!\Delta$ and $i_m\! =\! 0$ for some $m\! <\! k$, then $l_m\! =\! 0$ or $l_m\! =\! 1$ (but $l_m$ cannot be equal to $-1$).\smallskip

 Let $w\! =\! a_0\ldots a_{l-1}$ be a finite word over $\Si$. A sequence $r\! =\! (q_i,c_0^i,\ldots ,c_{k-1}^i)_{i<N}$ of configurations, where $N\! >\! l$, is called a {\bf run of $\mathcal{A}$ on $w$ starting in the configuration} 
$(p, c_0,\ldots ,c_{k-1})$ if
\begin{enumerate}
\ite[(1)]  $(q_0,c_0^0,\ldots ,c_{k-1}^0)\! =\! (p,c_0,\ldots, c_{k-1})$,
\ite[(2)] for each $i\! <\! N\! -\! 1$, there exists $b_i \in \Si \cup \{\lambda\}$ such that 
$$b_i: (q_i,c_0^i,\ldots ,c_{k-1}^i)\mapsto_{\mathcal{A}}(q_{i+1}, c_0^{i+1},\ldots ,c_{k-1}^{i+1})\mbox{,}$$
and $a_0\ldots a_{l-1}\! =\! b_0\ldots b_{N-2}$.
\end{enumerate}

 A run of $\mathcal{A}$ on $w$ starting in the configuration $(q_0,0,\ldots, 0)$ will be simply called a {\bf run of $\mathcal{A}$ on} $w$. The run is {\bf accepting} if it ends in a final state. The language $L(\mathcal{A})$ accepted by $\mathcal{A}$ is the set of finite words admitting an accepting run by $\mathcal{A}$.\smallskip
  
 Let $\sigma\! =\! a_0a_1\ldots$ be an $\om$-word over $\Si$. An $\om$-sequence of configurations 
$r\! =\! (q_i,c_0^i,\ldots ,c_{k-1}^i)_{i\in\omega}$ is called a {\bf run of $\mathcal{A}$ on $\sigma$ starting in the configuration} 
$(p,c_0,\ldots ,c_{k-1})$ if
\begin{enumerate}
\ite[(1)] $(q_0,c_0^0,\ldots ,c_{k-1}^0)=(p,c_0,\ldots, c_{k-1})$,

\ite[(2)] for each $i\!\in\!\omega$, there is $b_i\!\in\!\Si\cup\{\lambda\}$ such that 
$b_i\! :\! (q_i,c_0^i,\ldots ,c_{k-1}^i)\!\mapsto_{\mathcal{A}}\! (q_{i+1}, c_0^{i+1},\ldots ,c_{k-1}^{i+1})$, and either 
$b_0b_1\ldots\! =\! a_0a_1\ldots$, or $b_0b_1\ldots$ is a finite prefix of $a_0a_1\ldots$
\end{enumerate}

 The run $r$ is said to be {\bf complete} when $a_0a_1\ldots\! =\! b_0b_1\ldots$ For every such run, $\mbox{In}(r)$ is the set of all states entered infinitely often during the run $r$. A complete run of $\mathcal{A}$ on $\sigma$ starting in the configuration 
$(q_0,0,\ldots, 0)$ will be simply called a {\bf run of $\mathcal{A}$ on $\sigma$}. The \ol~accepted by $\mathcal{A}$ is 
\begin{center}
$L(\mathcal{A})\! :=\!\{\sigma\in\Si^\om\mid\mbox{ there exists a run }r~{ of }~\mathcal{A}\mbox{ on }\sigma\mbox{ such that } \mbox{In}(r)\cap F\!\neq\!\emptyset\}$.
\end{center}
\end{defi}

\noindent\bf Remark.\rm\ The acceptance condition for finite words here is by final states.  Some other acceptance conditions have been considered. In particular, in the last sections of the paper, we will consider acceptance of finite words by final states {\it and} counters having the value zero.


\section{$\!\!\!\!\!\!$ Topology}\indent

 We now recall some notions of topology, assuming  the reader to be familiar with the basic notions, that can be found in  \cite{Moschovakis80,Kechris94,Staiger97,PerrinPin}. The topological spaces in which we will work in this paper will be subspaces of $\Sigma^\omega$, where $\Si$ is either finite having at least two elements (like $2\! :=\!\{ {\bf 0},{\bf 1}\}$), or countably infinite. The topology on $\Sigma^\omega$ is the product topology of the discrete topology on $\Sigma$. For $w\!\in\!\Si^{<\om}$, the set 
$N_w\! :=\!\{\alpha\!\in\!\Sigma^\omega\mid w\!\subseteq\!\alpha\}$ is a basic clopen (i.e., closed and open) set of $\Sigma^\omega$. The open subsets of $\Sio$ are of the form $W^\frown\Si^\om\! :=\!\{ w\sigma\mid w\!\in\! W\mbox{ and }\sigma\!\in\!\Si^\om\}$, where $W\!\subseteq\!\Si^{<\om}$. When $\Si$ is finite, this topology is called the {\bf Cantor topology} and $\Sio$ is compact. When $\Si\! =\!\omega$, 
$\Sigma^\omega$ is the Baire space, which is homeomorphic to 
$\mathbb{P}_\infty\! :=\!\{\alpha\!\in\! 2^\omega\mid\forall i\!\in\!\omega\ \exists j\!\geq\! i\ \ \alpha (j)\! =\! {\bf 1}\}$, via the map defined on $\omega^\omega$ by $h(\beta )\! :=\! {\bf 0}^{\beta (0)}{\bf 1}{\bf 0}^{\beta (1)}{\bf 1}\ldots$ There is a natural metric on $\Sio$, the {\bf prefix metric} defined as follows. For $\sigma\!\not=\!\tau\!\in\!\Sio$, 
$d(\sigma ,\tau )\! :=\! 2^{-l_{pref(\sigma ,\tau )}}$, where $l_{pref(\sigma ,\tau )}$ is the first natural number $n$ such that 
$\sigma (n)\!\not=\!\tau (n)$. The topology induced on $\Sigma^\omega$ by this metric is our topology.\medskip

 We now define the Borel hierarchy.

\begin{defi}
Let $X$ be a topological space, and $n\!\geq\! 1$ be a natural number. The classes ${\bf\Si}_n^0(X)$ and ${\bf\Pi}_n^0(X) $ of the {\bf Borel hierarchy} are inductively defined as follows:\smallskip

${\bf \Si}^0_1(X) $ is the class of open subsets of $X$.\smallskip

${\bf \Pi}^0_1(X) $ is the class of closed subsets of $X$.\smallskip

${\bf \Si}^0_{n+1}(X)$ is the class of countable unions of ${\bf \Pi}^0_n$-subsets of  $X$.\smallskip

${\bf \Pi}^0_{n+1}(X)$ is the class of countable intersections of ${\bf \Si}^0_n$-subsets of $X$.

\noindent The Borel hierarchy is also defined for the transfinite levels. Let $\xi\!\geq\! 2$ be a countable ordinal.\smallskip

${\bf \Si}^0_\xi (X)$ is the class of countable unions of subsets of $X$ in $\bigcup_{\gamma <\xi}~{\bf \Pi}^0_\gamma$.\smallskip

${\bf \Pi}^0_\xi (X)$ is the class of countable intersections of subsets of $X$ in $\bigcup_{\gamma <\xi}~{\bf \Si}^0_\gamma$.
\end{defi}

 Suppose now that $\xi\!\geq\! 1$ is a countable ordinal and $X\!\subseteq\! Y$, where $X$ is equipped with the induced topology. Then 
$\boraxi (X)\! =\!\{ A\cap X\mid A\!\in\!\boraxi (Y)\}$, and similarly for $\bormxi$, see \cite[Section 22.A]{Kechris94}. Note that we defined the Borel classes ${\bf\Si}^0_\xi(X)$ and ${\bf\Pi}^0_\xi (X)$ mentioning the space $X$. However, when the context is clear, we will sometimes omit $X$ and denote ${\bf \Si}^0_\xi(X)$ by ${\bf \Si}^0_\xi$ and similarly for the dual class. The Borel classes are closed under finite intersections and unions, and continuous preimages. Moreover, $\boraxi$ is closed under countable unions, and $\bormxi$ under countable intersections. As usual, the ambiguous class $\borxi$ is the class $\boraxi\cap\bormxi$. The class of {\bf Borel sets} is $\borel\! :=\!\bigcup_{1\leq\xi <\omega_1}\ \boraxi\! =\!\bigcup_{1\leq\xi <\omega_1}\ \bormxi$, where $\om_1$ is the first uncountable ordinal. The {\bf Borel hierarchy} is as follows:
$$\begin{array}{ll}  
& \ \ \ \ \ \ \ \ \ \ \ \ \ \ \ \ \ \ \ \ \ \ \ \ \ \boraone\! =\!\hbox{\rm open}\ \ \ \ \ \ \ \ \ \ \ \ \ 
\boratwo\! \ \ \ \ \ \ \ \  \ \ \ 
\ldots\ \ \ \ \ \ \ \ \ \ \ \ \boraom\ \ \ \ \ \ldots\cr  
& \borone\! =\!\hbox{\rm clopen}\ \ \ \ \ \ \ \ \ \ \ \ \ \ \ \ \ \ \ \ \ \ \ \ \ \ \ 
\bortwo\ \ \ \ \ \ \ \ \ \ \ \ \ \ \ \ \ \ \ \ \ \ \ \ \ \ \ \ \ \ \ \borom\ \ \ \ \ \ \ \ \ \ \ \ \ \ \ \ \ \ \ \ \ \ \borel\cr
& \ \ \ \ \ \ \ \ \ \ \ \ \ \ \ \ \ \ \ \ \ \ \ \ \ \bormone\! =\!\hbox{\rm closed}\ \ \ \ \ \ \ \ \ \ \bormtwo\! \ \ \ \  \ \ \ \ \ \ \ \ \ldots
\ \ \ \ \ \ \ \ \ \ \ \ \bormom\ \ \ \ \ \ldots
\end{array}$$
This picture means that any class is contained in every class to the right of it, 
and the inclusion is strict in any of the spaces $\Sigma^\omega$. A subset of $\Si^\om$ is a Borel set of {\bf rank} $\xi$ if 
it is in ${\bf\Si}^0_\xi\cup {\bf\Pi}^0_\xi$ but not in $\bigcup_{1\leq\gamma <\xi}~({\bf\Si}^0_\gamma\cup {\bf\Pi}^0_\gamma)$.\medskip

 We now define completeness with respect to reducibility by continuous functions. Let ${\bf\Gamma}$ be a class of sets of the form $\boraxi$ or $\bormxi$. A subset $C$ of $\Si^\om$ is said to be ${\bf\Gamma}$-{\bf complete} if $C$ is in 
${\bf\Gamma}(\Si^\om )$ and, for any finite alphabet $Y$ and any $A\!\subseteq\! Y^\om$, $A\!\in\! {\bf\Gamma}$ if and only if  there exists a continuous function $f\! :\! Y^\om\!\ra\!\Si^\om$ such that $A\! =\! f^{-1}(C)$. The ${\bf\Si}^0_n$-complete sets and the ${\bf\Pi}^0_n$-complete sets are thoroughly characterized in \cite{Staiger86a}. Recall that a subset of $\Sio$ is ${\bf\Si}^0_\xi$ (respectively ${\bf\Pi}^0_\xi$)-complete if it is in ${\bf\Si}^0_\xi$ but not in ${\bf\Pi^0_\xi}$ (respectively in ${\bf\Pi}^0_\xi$ but not in 
${\bf\Si}^0_\xi$), \cite{Kechris94}. For example, the singletons of $2^\omega$ are $\bormone$-complete. The set 
$\mathbb{P}_\infty$ defined at the beginning of the present section is a well known example of a $\bormtwo$-complete set.\medskip

 The class $\check {\bf\Gamma}\! :=\!\{\neg A\mid A\!\in\! {\bf\Gamma}\}$ is the class of the complements of the sets in ${\bf\Gamma}$. In particular, $\check {{\bf\Si}^0_\xi}\! =\! {\bf\Pi}^0_\xi$ and 
$\check {{\bf\Pi}^0_\xi}\!=\! {\bf\Si}^0_\xi$.\medskip

 There are some subsets of the topological space $\Sio$ which are not Borel sets. In particular, there is another hierarchy beyond the Borel hierarchy,  called the projective hierarchy. The first class of the projective hierarchy is the class ${\bf\Si}^1_1$ of analytic sets. A subset $A$ of $\Sio$ is {\bf analytic} if we can find a finite alphabet $Y$ and a Borel subset $B$ of 
$(\Si\!\times\! Y)^\om$ such that $x\!\in\! A\Leftrightarrow\exists y\!\in\! Y^\om ~(x, y)\!\in\! B$, where 
$(x, y)\!\in\! (\Si\!\times\! Y)^\om$ means that $(x, y)(i)\! =\!\big( x(i),y(i)\big)$ for each natural number $i$.

\vfill\eject

 A subset of $\Sigma^\omega$ is analytic if it is empty, or the image of the Baire space by a continuous map. The class ${\bf\Si}^1_1$ of analytic sets contains the class of Borel sets in any of the spaces $\Sigma^\omega$. Note that ${\bf\Delta}_1^1\! =\! {\bf\Si}^1_1 \cap {\bf\Pi}^1_1$, where ${\bf\Pi}^1_1\! :=\!\check\ana$ is the class of {\bf co-analytic} sets, i.e., of complements of analytic sets.\medskip
 
 The $\om$-power of a finitary language $L$ is always an analytic set. Indeed, if $L$ is finite and has $n$ elements, then $L^\om$  is the continuous image of the compact set 
$\{ {\bf 0},{\bf 1},\ldots ,{\bf n\! -\! 1}\}^\om$. If $L$ is infinite, then there is a bijection between $L$ and $\om$, and $L^\om$ is the continuous image of the Baire space $\om^\om$, \cite{Simonnet92}.

\section{$\!\!\!\!\!\!$ The operation  ``exponentiation of sets"}\indent

 The Wadge hierarchy of Borel sets is a great refinement of the Borel hierarchy. Wadge gave first a description of this hierarchy, see \cite{Wadge83}. Duparc  got in  \cite{Duparc01} a new proof of Wadge's results for the case of  Borel sets of finite rank, and he gave a normal form of  Borel sets of finite rank, i.e., an inductive construction of a Borel set of every given degree. His proof relies on set theoretic operations which are the counterpart of arithmetical operations over ordinals needed to compute the Wadge degrees.\medskip
 
 In fact J. Duparc studied the Wadge hierarchy via the study of the conciliating hierarchy. He introduced in \cite{Duparc01} the conciliating sets, which are sets of finite {\it or} infinite words over an alphabet $\Si$, i.e. subsets of $\Si^{<\om}\cup\Si^\om\! =\!\Si^{\leq\om}$. Among the  set theoretic operations which are  defined over concilating sets, we shall only 
need in this paper the operation of exponentiation. We first recall the following.

\begin{defi}
Let  $\Si_A$ be a finite alphabet, $\tla$ be a letter out of $\Si_A$, $\Si\! :=\!\Si_A\cup\{\tla\}$, and $x$ be a finite or infinite word over the alphabet $\Si$. Then $x^\tla$ is inductively defined as follows.\smallskip

- $\lambda^\tla\! :=\!\lambda$.\smallskip

- For a finite word $u\!\in\!\Si^{<\om}$, 
$\left\{\!\!\!\!\!\!\!
\begin{array}{ll}
& (ua)^\tla\! :=\! u^\tla a\mbox{ if }a\!\in\!\Si_A\mbox{,}\cr
& (u\tla)^\tla\! :=\! u^\tla\mbox{ with its last letter removed if }|u^\tla|\! >\! 0\mbox{,}\cr
& (u\tla)^\tla\! :=\!\lambda\mbox{ if }|u^\tla|\! =\! 0.
\end{array}
\right.$\smallskip

- For an infinite word $\sigma$, $\sigma^\tla\! :=\!\mbox{lim}_{n\in\om}~(\sigma\vert n)^\tla$, where, given 
$(w_n)\!\in\! (\Si_A^{<\om})^\om$ and $w\!\in\!\Si_A^{<\om}$, 
$$w\!\subseteq\!\mbox{lim}_{n\in\om}~w_n\Leftrightarrow\exists p\!\in\!\omega ~~\fa n\!\geq\! p~~w_n\vert |w|\! =\! w.$$
\end{defi}

\noindent\bf Remark.\rm\ For $x\!\in\!\Si^{\leq\om}$, $x^\tla$ denotes the string $x$, once every $^\tla$ occuring in $x$ has been ``evaluated" as the back space operation (the one familiar to your computer!), proceeding from left to right inside $x$. In other words, $x^\tla\! =\! x$ from which every interval of the form $``a\tla "$ ($a\!\in\!\Si_A$) is removed.\medskip

 For example, if $x\! =\! (a\tla)^n$ for some $n\!\geq\! 1$, $x\! =\! (a\tla)^\om$ or $x\! =\! (a\tla\tla)^\om$ then $x^\tla\! =\!\lambda$. If  $x\! =\! (ab\tla)^\om$, then $x^\tla\! =\! a^\om$. If $x\! =\! bb(\tla a)^\om$, then $x^\tla\! =\! b$.\medskip 

 We now can define the operation $A\!\mapsto\! A^\sim$ of exponentiation of conciliating sets.

\begin{defi}
Let  $\Si_A$ be a finite alphabet, $\tla$ be a letter out of $\Si_A$, $\Si\! :=\!\Si_A\cup\{\tla\}$, and $A\!\subseteq\!\Si_A^{\leq\om}$. Then we set $A^\sim\! :=\!\{ x\!\in\!\Si^{\leq\om}\mid x^\tla\!\in\! A\}$.
\end{defi}

 The operation $\sim$ is monotone with regard to the Wadge ordering and produces some sets of higher complexity. Duparc considered the following correspondence. If $\Si_A$ is a finite alphabet, $A\!\subseteq\!\Si_A^{\leq \om}$ and $d$ is a letter out of $\Si_A$, then we define
$$A^d\! :=\!\{\sigma\!\in\! (\Si_A\cup\{ d\} )^\om\mid\sigma (/d)\!\in\! A\}\mbox{,}$$
where $\sigma (/d)$ is the sequence obtained from $\sigma$ by removing every occurrence of the letter $d$.\medskip

 We recall the results useful in this paper. 

\begin{thm} [Duparc \cite{Duparc01} ]\label{thedup} Let  $\Si_A$ be a finite alphabet.
\begin{itemize}
\ite[(a)] 
Let $A\!\subseteq\!\Si_A^{\leq \om}$, and $n\!\geq\! 1$ be a natural number. If $A^d\!\subseteq\! (\Si_A\cup\{ d\} )^\om$ is 
${\bf\Si}_n^0$-complete (respectively, ${\bf \Pi}_n^0$-complete), then $(A^\sim)^d$ is ${\bf \Si}_{n+1}^0$-complete (respectively, 
${\bf \Pi}_{n+1}^0$-complete).
\ite[(b)] Let $A\!\subseteq\!\Si_A^\om$, and $n\!\geq\! 2$ be a natural number. If $A$ is ${\bf \Pi}_n^0$-complete, then $A^\sim$ is 
${\bf \Pi}_{n+1}^0$-complete.
\end{itemize}
\end{thm}

\noindent\bf Remark.\rm\ Item (b) of the  preceding theorem follows from (a) because\medskip

- whenever $A\!\subseteq\!\Si_A^\om$, $n\!\geq\! 2$ is a natural number and $A$ is ${\bf \Pi}_n^0$-complete, then $A^d$ is also 
${\bf \Pi}_n^0$-complete,\smallskip

- whenever $A\!\subseteq\!\Si_A^\om$, $n\!\geq\! 3$ is a natural number and $A^d\!\subseteq\! (\Si_A\cup\{ d\} )^\om$ is 
${\bf \Pi}_n^0$-complete, then $A$ is also a ${\bf \Pi}_n^0$-complete set.\medskip

 This property was useful in \cite{Fin01a} to study the $\om$-powers of finitary languages. The first author proved in \cite{Fin01a} that the class $CFL_\om$ of context-free $\om$-languages, (i.e., those which are accepted by  B\"uchi pushdown automata), is closed under this operation $\sim$.\medskip 
  
 We now recall a  slightly modified variant of the operation $\sim$, introduced in \cite{Fin01a}, and which is particularly suitable to infer properties of $\om$-powers. 
  
\begin{defi}
Let  $\Si_A$ be a finite alphabet, $\tla$ be a letter out of $\Si_A$, $\Si\! :=\!\Si_A\cup\{\tla\}$, and $A\!\subseteq\!\Si_A^{\leq\om}$. Then we set $A^\approx\! :=\!\{ x\!\in\!\Si^{\leq\om}\mid x^\tla\!\in\! A\}$, where $x^\tla$ is inductively defined as follows.\smallskip

- $\lambda^\tla\! :=\!\lambda$.\smallskip

- For a finite word $u\!\in\!\Si^{<\om}$, 
$\left\{\!\!\!\!\!\!\!
\begin{array}{ll}
& (ua)^\tla\! :=\! u^\tla a\mbox{ if }a\!\in\!\Si_A\mbox{,}\cr
& (u\tla)^\tla\! :=\! u^\tla\mbox{ with its last letter removed if }|u^\tla|\! >\! 0\mbox{,}\cr
& (u\tla)^\tla\mbox{ is undefined if }|u^\tla|\! =\! 0.
\end{array}
\right.$\smallskip

- For an infinite word $\sigma$, $\sigma^\tla\! :=\!\mbox{lim}_{n\in\om}~(\sigma\vert n)^\tla$.
\end{defi}

 The only difference is that here $(u\tla)^\tla$ is undefined if $|u^\tla|\! =\! 0$. It is easy to see that if $A\!\subseteq\!\Si_A^\om$ is a Borel set such that $A\!\neq\!\Si_A^\om$, i.e., $A^-\!\neq\!\emptyset$, 
then $A^\approx$ is Wadge equivalent to $A^\sim$ (see \cite{Fin01a}) and thus one can get the following version of Theorem \ref{thedup}.(b). 

\begin{thm}\label{thedup2}
Let  $\Si_A$ be a finite alphabet, and $n\!\geq\! 2$ be a natural number. If $A\!\subseteq\!\Si_A^\om$ is ${\bf \Pi}_n^0$-complete, 
then $A^\approx$ is ${\bf \Pi}_{n+1}^0$-complete.
\end{thm}

 

  
\section{$\!\!\!\!\!\!$ ${\bf \Pi}_n^0$-complete $\om$-powers}
  
\bf Notation.\rm\ Let $\Si_A$ be a finite alphabet, $\tla$ be a letter out of $\Si_A$, and $\Si\! :=\!\Si_A\cup\{\tla\}$. The language $L_3$ over $\Si$ is the context-free language generated by the context free grammar with the following production rules:
$$\begin{array}{ll}
& S\!\ra\! aS\tla S\mbox{ with }a\!\in\!\Si_A\mbox{,}\cr
& S\!\ra\! a\tla S\mbox{ with }a\!\in\!\Si_A\mbox{,}\cr
& S\!\ra\!\lambda 
\end{array}$$
(see [H-U] for the basic notions about grammars). This language $L_3$ corresponds to the words where every letter of $\Si_A$ has been removed after using the backspace operation. It is easy to see that $L_3$ is a \de one-counter \la , i.e., $L_3$ is accepted by a \de one-counter automaton. Moreover,  for $a\!\in\!\Si_A$, the language $L_3a$ is also accepted by a deterministic one-counter automaton.\medskip

 We can now state the following result. 
 
\begin{lem}[see \cite{Fin01a}] Whenever $A\!\subseteq\!\Si_A^\om$, the \ol~ $A^\approx\!\subseteq\!\Si^\om$ is obtained by substituting in $A$ the language $L_3a$ for each letter $a\!\in\!\Si_A$.\end{lem} 
  
   An $\om$-word $\sigma\!\in\! A^\approx$ may be considered as an $\om$-word $\sigma^\tla\!\in\! A$ to which we possibly add, before the first letter $\sigma^\tla (0)$ of $\sigma^\tla$ (respectively, between two consecutive letters $\sigma^\tla (n)$ and 
$\sigma^\tla (n\! +\! 1)$ of $\sigma^\tla$), a finite word belonging to the context free (finitary) language $L_3$. 
  
\begin{cor} \label{presop}
Whenever $A\!\subseteq\!\Si_A^\om$ is an $\om$-power of a \la $L_A$, i.e., $A\! =\! L_A^\infty$, then $A^\approx$ is also an 
$\om$-power, i.e., there exists a (finitary) language $E_A$ such that $A^\approx\! =\! E_A^\infty$. Moreover, if the language $L_A$ is in the class $OCL(k)$ for some natural number $k$, then the language $E_A$ can be found in the class $OCL(k\! +\! 1)$. 
\end{cor}

\proo  Let $h\! :\!\Si_A\!\rightarrow\! 2^{\Si^{<\omega}}$ be the substitution defined by $a\!\mapsto\! L_3a$, where $L_3$ is the context free  language defined above. Then it is easy to see that  now $A^\approx$ is obtained by substituting in $A$ the language
$L_3a$ for each letter $a\!\in\!\Si_A$. Thus $E_A\! =\! h(L_A)$ satisfies the statement of the theorem.\hfill{$\square$}\medskip

 We now  recall the following result, proved in \cite{Fin01a}. 

\begin{thm}\label{pin} 
For each natural number $n\!\geq\! 1$, there is a context free language $P_n$ in the subclass of iterated counter languages such that $P_n^\infty$ is ${\bf \Pi}_n^0$-complete.
\end{thm}

\proo
 Let $B_1\! =\!\{\sigma\!\in\!\{ {\bf 0},{\bf 1}\}^\om\mid\fa i\!\in\!\om ~~\sigma (i)\! =\! {\bf 0}\}\! =\!\{ {\bf 0}\}^\infty$. $B_1$ is a ${\bf \Pi}_1^0$-complete set of the form $P_1^\infty$ where $P_1$ is the singleton containing only the word ${\bf 0}$. Note that that $P_1\! =\!\{ {\bf 0}\}\! =:\! 1$ is a regular  language, hence in the class $OCL(0)$.\medskip

 Let then $B_2\! =\!\mathbb{P}_\infty$ be the well known ${\bf \Pi}_2^0$-complete regular $\om$-language. Note that 
$B_2\! =\! (1^{<\om}{\bf 1})^\infty$. Let $P_2\! :=\! 1^{<\om}{\bf 1}$. Then $P_2$ is  a regular  language, hence in the class $OCL(0)$.\medskip
 
 We now consider the substitution $h\! :\!\{ {\bf 0},{\bf 1}\}\!\ra\! 2^{(\{ {\bf 0},{\bf 1}\}\cup\{\tla\})^{<\omega}}$ from the proof of Corollary \ref{presop}, and set 
$P_3\! :=\! h(P_2)$, which is a context-free language in the class  $OCL(1)$. Note that the set 
${P_3^\infty\! =\! h(P_2)^\infty\! =\! (P_2^\infty )^\approx}$ is ${\bf \Pi}_3^0$-complete, by Theorem \ref{thedup}.\medskip

 Iterating this method $n\!\geq\! 1$ times, we easily obtain a context free language $P_{n+2}\!\in\! OCL(n)$ such 
that $P_{n+2}^\infty$ is ${\bf \Pi}_{n+2}^0$-complete.\hfill{$\square$}

\vfill\eject

 Note that $P_1$ and $P_2$ are regular, hence accepted by some (real-time deterministic) finite automata (without any counter). On the other hand, the languages $L_3a$, for $a\!\in\!\Si_A$, are one-counter languages. Moreover we can easily see that, for each $a\!\in\!\Si_A$,  the language $L_3a$ is accepted by a real-time one-counter automaton $\mathcal{A}$, such that 
$\mathcal{A}$ accepts a finite word $w$ iff there is a run on $w$ ending in an accepting state {\it and} with empty stack. This implies that the language $P_3$ is also accepted by a real-time one-counter automaton, by final states {\it and} empty stack. We can now state the following proposition. 

\begin{prop} \label{sub} 
Let $A\!\subseteq\!\Si_A^{<\om}$ be a finitary language accepted by a real-time one-counter automaton $\mathcal{A}$ accepting words by final states {\it and} empty stack, and let $h\! :\!\Si_A\!\ra\! 2^{\Si^{<\om}}$ be the substitution defined by 
$a\!\mapsto\! L_3a$, where $L_3$ is the one-counter language defined above. Then the language $h(A)$ is in $OCL(2)$, and it is accepted by a real-time two-iterated counter automaton $\mathcal{B}$ accepting words by final states {\it and} empty stack. 
\end{prop}
    
 We explain in an informal way the idea of the construction of the real-time two-iterated counter automaton $\mathcal{B}$ from the automaton $\mathcal{A}$. The stack alphabet of $\mathcal{A}$ is of the form $\{Z_0, z_0\}$, and the stack alphabet of 
$\mathcal{B}$ is of the form $\{Z_0, z_0, z_1\}$. The automaton $\mathcal{B}$ starts the reading of a word over the alphabet $\Si$     as the one-counter automaton $\mathcal{A}$ accepting the language $A$. Then at any moment of the computation it may guess (using the non-determinism) that it reads a finite segment $w$ of $L_3$ that will be erased (using the eraser $\tla$). It reads $w$ using the additional stack letter $z_1$ which permits to simulate a one-counter automaton at the top of the stack while keeping the memory of the stack of $\mathcal{A}$ (which is actually a counter). Then, after the reading of $w$, $\mathcal{B}$ simulates again the one-counter automaton $\mathcal{A}$, and so on. The automaton 
$\mathcal{B}$ accepts words by final states (corresponding to final states of $\mathcal{A}$) {\it and} empty stack. We now state one of our main technical results. 

\begin{prop} \label{21} 
Let $A\!\subseteq\!\Si^{<\om}$ be a finitary language accepted by a real-time two-iterated counter automaton $\mathcal{A}$ accepting words by final states  {\it and} empty stack, and such that the $\om$-power $A^\infty$ is ${\bf \Si}_n^0$-complete  (respectively, ${\bf \Pi}_n^0$-complete) for some natural number $n\!\geq\! 3$. Then we can find a finite alphabet $Y$ and a finitary language $B\!\subseteq\! Y^{<\om}$ such that $B$ is accepted by a real-time one-counter automaton $\mathcal{B}$ accepting words by final states {\it and} empty stack, and $B^\infty$ is ${\bf\Si}_n^0$-complete (respectively,  ${\bf\Pi}_n^0$-complete). Moreover one can take a two-letter alphabet $Y=\{0, 1\}$ with the same property. 
\end{prop}

\proo Note first that we have already seen in Section \ref{section-automata} that a (real-time) two-iterated counter automaton (accepting words by final states  {\it and} empty stack) may be seen as a (real-time) two-counter automaton (accepting words by final states {\it and} counters having the value zero). The idea is to code the content of two counters. We shall need the following notion. Let $m\!\geq\! 1$ be a natural number, and $n,p,q$ be natural numbers such that neither $2$ nor $3$ divides $n\!\geq\! 1$, and $m\! =\! n.2^p.3^q$. Then we set $M_2(m)\! :=\! p$ and $M_3(m)\! =\! q$. So $2^{M_2(m)}$ is the greatest power of $2$ which divides $m$, and $2^{M_3(m)}$ is the greatest power of $3$ which divides $m$.\medskip 

Let then $\mathcal{A}\! :=\! (Q,\Si,q_0,\Delta ,F)$ be a (real-time) two-iterated counter automaton accepting the language 
$A\! =\! L(\mathcal{A})\!\subseteq\!\Sis$, by final states  {\it and} empty stack. We define the finitary language $\mathcal{L}$ as the set of finite words over the alphabet $\Si\cup\{ {\bf 0},{\bf 1},{\bf 2}\}$, where  ${\bf 0},{\bf 1},{\bf 2}$ are new letters not in $\Si$, of the form 
${{}^\frown}_{i<n}~v_i~a_i~{\bf 1}~w_i~z_i~{\bf 2}~u_{i+1}$, where $|v_0|\! =\! 1$, $n\!\geq\! 1$, $v_i, w_i\!\in\! 1^+$, 
$z_i,u_i\!\in\! 1^{<\om}$, $a_i\!\in\!\Si$, $|u_{i+1}|\! =\! |z_i|$ and we can find a sequence $(q_i)_{i\leq n}$ of states of $Q$ and integers $l_i,l'_i\!\in\!\{ -1,0,1\}$ such that, for each $i\! <\! n$,
$$a_i\! :\!\big( q_i,M_2(|v_i|),M_3(|v_i|)\big)\!\mapsto_{\mathcal{A}}\!\big( q_{i+1},M_2(|v_i|)\! +\! l_i,M_3(|v_i|)\! +\! l'_i\big)$$
and $|w_i|\! =\! |v_i|.2^{l_i}.3^{l'_i}$.

\vfill\eject

 Moreover, the state $q_n$ is a final state of $\mathcal{A}$, i.e., $q_n\!\in\! F$, and 
$M_2(|w_{n-1}|)\! =\! 0$, $M_3(|w_{n-1}|)\! =\! 0$. Note that the state $q_0$ of the sequence $(q_i)_{i\leq n}$ is also the initial state of $\mathcal{A}$.\medskip 

\noindent\bf Claim 1\label{lem4-}\it\ 
The language $\mathcal{L}$ is accepted by a one-counter automaton $\mathcal{C}$.\rm\medskip

 Indeed, we shall explain informally the behaviour of a one-counter automaton $\mathcal{C}$ accepting the finitary language $\mathcal{L}$.  We first consider the reading of a word $w\!\in\! (\Si\cup\{ {\bf 0},{\bf 1},{\bf 2}\})^{<\om}$ of the form 
$$({{}^\frown}_{i<n}~ {\bf 0}^{p_i}~a_i~{\bf 1}~ {\bf 0}^{m_i}~{\bf 2})~{\bf 0}^{p_n}\mbox{,}$$
where the $p_i,m_i,$'s are positive integers, and the $a_i$'s are in $\Si$.\medskip

 Using the finite control, (i.e.,  a finite set of states and a set of transitions involving only these states, and the input letters that are read; this corresponds to the behaviour of a finite state automaton)   the automaton $\mathcal{C}$ first checks that the first three letters of $w$ form an initial segment ${\bf 0}~a_0~{\bf 1}$, for some letter $a_0\!\in\!\Si$. Moreover, when reading the $p_0\! =\! 1$ letter ${\bf 0}$ before $a_0$, the automaton $\mathcal{C}$, using the finite control, checks that $p_0\! >\! 0$ and determines whether 
$M_2(p_0)\! =\! 0$, and whether $M_3(p_0)\! =\! 0$. Here we actually have $M_2(p_0)\! =\! 0$ and $M_3(p_0)\! =\! 0$. Moreover the counter content of $\mathcal{C}$ is increased by one for each letter ${\bf 0}$ read.\medskip
  
  The automaton $\mathcal{C}$ now reads the letter $a_0$ and it guesses a transition of $\mathcal{A}$ leading to  
$$a_0\! :\!\big( q_0,M_2(p_0),M_3(p_0)\big)\!\mapsto_{\mathcal{A}}\!\big( q_1,M_2(p_0)\! +\! l_0,M_3(p_0)\! +\! l'_0\big)$$
We set $v_0\! :=\! {\bf 0}^{p_0}$ and $w_0\! :=\! {\bf 0}^{p_0.2^{l_0}.3^{l'_0}}$.\medskip
 
 The counter value is now equal to $p_0$ and, when reading the letters ${\bf 0}$ following $a_0$, the automaton $\mathcal{C}$ checks that $m_0\!\geq\! p_0.2^{l_0}.3^{l'_0}$ in such a way that the counter value becomes zero after having read the 
$p_0.2^{l_0}.3^{l'_0}$ letters ${\bf 0}$ following the first letter ${\bf 1}$. For instance, if $l_0\! =\! l'_0\! =\! 1$, then 
$|w_0|\! =\! |v_0|.6$ so this can be done by decreasing the counter content by one each time six letters ${\bf 0}$ are read. The other cases are treated similarly. The details are here left to the reader.\medskip

Note also that the automaton $\mathcal{C}$ has kept in its finite control the value of the state $q_1$.\medskip

 We now set ${\bf 0}^{m_0}\! :=\! w_0.z_0$. We have seen that, after having read $w_0$, the counter value of the automaton 
$\mathcal{C}$ is equal to zero. Now when reading $z_0$ the counter content is increased by one for each letter read so that it becomes $|z_0|$ after having read $z_0$. The automaton $\mathcal{C}$ now reads a letter $A$ and next decreases its counter by one for each letter ${\bf 0}$ read until the counter content is equal to zero. We set ${\bf 0}^{p_1}\! :=\! u_1.v_1$ with 
$u_1\! =\! z_0$. The automaton $\mathcal{C}$ now reads the segment $v_1$. Using the finite control, it checks that 
$|v_1|\! >\! 0$ and determines whether $M_2(|v_1|)\! =\! 0$, and whether $M_3(|v_1|)\! =\! 0$. Moreover the counter content is increased by one for each letter ${\bf 0}$ read. The automaton $\mathcal{C}$ now reads the letter $a_1$ and it guesses a transition of $\mathcal{A}$ leading to  
$$a_1\! :\!\big( q_1,M_2(|v_1|),M_3(|v_1|)\big)\!\mapsto_{\mathcal{A}}
\!\big( q_2,M_2(|v_1|)\! +\! l_1,M_3(|v_1|)\! +\! l'_1\big)$$
We set $w_1\! :=\! {\bf 0}^{|v_1|.2^{l_1}.3^{l'_1}}$. The counter value is now equal to $|v_1|$. The automaton $\mathcal{C}$ now reads the second letter ${\bf 1}$ and, when reading the $m_1$ letters ${\bf 0}$ following this letter ${\bf 1}$, the automaton 
$\mathcal{C}$ checks that $m_1\!\geq\! |v_1|.2^{l_1}.3^{l'_1}$ in such a way that the counter value becomes zero after having read the $|v_1|.2^{l_1}.3^{l'_1}$ letters ${\bf 0}$ following the second letter ${\bf 1}$.

\vfill\eject

 For instance, if $l_1\! =\! 0$ and $l'_1\! =\! -1$, then $|w_1|\! =\! |v_1|.3^{-1}$, so this can be done by decreasing the counter content by three each time one letter ${\bf 0}$ is read. And if $l_1\! =\! -1$ and $l'_1\! =\! -1$, then $|w_1|\! =\! |v_1|.2^{-1}.3^{-1}\! =\! |v_1|.6^{-1}$ so this can be done by decreasing the counter content by six each time one letter ${\bf 0}$ is read. The other cases are treated similarly. The details are here left to the reader.\medskip

 Note that these different cases can be treated using $\lambda$-transitions, in such a way that there will be at most 5 consecutive 
$\lambda$-transitions during a run of $\mathcal{C}$ on $w$. This will be an important useful fact in the sequel.\medskip
 
 Note also that the automaton $\mathcal{C}$ has kept in its finite control the value of the state $q_2$.\medskip
 
 The reading of $w$ by $\mathcal{C}$ continues similarly. An acceptance condition by final states and empty stack  can be used to ensure that  $q_n\!\in\! F$, $M_2(|w_{n-1}|)\! =\! 0$, $M_3(|w_{n-1}|)\! =\! 0$, and $|z_{n-1}|\! =\! p_n$.\medskip 
 
 In order to complete the proof, we can remark that 
$\mathcal{R}\! =\! ({\bf 0}^{<\om}.\Si.{\bf 1}.{\bf 0}^{<\om}.{\bf 2}.{\bf 0}^{<\om})^{<\om}$ is a regular language,  so we have only considered the reading by $\mathcal{C}$ of words $w\!\in\!\mathcal{R}$. Indeed, if the language $L(\mathcal{C})$ was not included into 
$\mathcal{R}$, then we could replace it with $L(\mathcal{C})\cap\mathcal{R}$ because the class $OCL$ is closed under intersection with regular languages (by a classical construction of product of automata, the language $\mathcal{R}$ being accepted by a deterministic finite-state automaton).\hfill{$\diamond$}\medskip  

 We will now use a coding of $\om$-words over $\Si$ given by the map $g_{N,l}\! :\!\Sio\!\ra\! (\Si\cup\{ {\bf 0},{\bf 1},{\bf 2}\})^\om$, where $(N,l)\!\in\!\mathcal{P}\! :=\!\{ (N,l)\!\in\! (\om\!\setminus\!\{ 0\} )\!\times\!\om\mid 6\mbox{ does not divide }N\}$, and 
$$g_{N,l}(\sigma )\! :=\! {\bf 0}~({{}^\frown}_{i\in\om}~\sigma (i)~{\bf 1}~{\bf 0}^{N.6^{l+i}}~{\bf 2}~{\bf 0}^{N.6^{l+i}}).$$
\bf Claim 2\label{lem5}\it\ The equality $\big( L(\mathcal{A})\big)^\infty\! =\! g_{N,l}^{-1}(\mathcal{L}^\infty )$ holds, i.e., 
$\fa\sigma\!\in\!\Si^{\om}$, $g_{N,l}(\sigma )\!\in\!\mathcal{L}^\infty\Leftrightarrow\sigma\!\in\!\big( L(\mathcal{A})\big)^\infty$.\rm\medskip

 Indeed, let $\mathcal{A}$ be a real-time two-iterated counter automaton accepting finite words over $\Si$, by final states 
{\it and} empty stack, and $\mathcal{L}\!\subseteq\! (\Si\cup\{ {\bf 0},{\bf 1},{\bf 2}\})^{<\om}$ be defined as above.\medskip

 Let $\sigma\!\in\!\Si^{\om}$ be an $\om$-word such that $g_{N,l}(\sigma )\!\in\!\mathcal{L}^\infty$. Recall that $g_{N,l}(\sigma )$ can be written 
$${\bf 0}~({{}^\frown}_{i\in\om}~\sigma (i)~{\bf 1}~{\bf 0}^{N.6^{l+i}}~{\bf 2}~{\bf 0}^{N.6^{l+i}})$$
As $g_{N,l}(\sigma )\!\in\!\mathcal{L}^\infty$, we can also write 
$g_{N,l}(\sigma )\! =\! {^\frown}_{j\in\om}~({{}^\frown}_{i<n_j}~v_{i,j}~a_{i,j}~{\bf 1}~w_{i,j}~z_{i,j}~{\bf 2}~u_{i+1,j})$, where, for each natural number $j$, $|v_{0,j}|\! =\! 1$, $n_j\!\geq\! 1$, $v_{i,j},w_{i,j}\!\in\! 1^+$, $z_{i,j},u_{i,j}\!\in\! 1^{<\om}$, 
$a_{i,j}\!\in\!\Si$, $|u_{i+1,j}|\! =\! |z_{i,j}|$ and we can find a sequence $(q_{i,j})_{i\leq n_j}$ of states of $Q$ such that 
$q_{0,j}\! =\! q_0$ is the initial state of $\mathcal{A}$, and integers $l_{i,j}, l'_{i,j}\!\in\!\{-1,0,1\}$ such that, for each $i\! <\! n_j$, 
$$a_{i,j}\! :\!\big( q_{i,j},M_2(|v_{i,j}|),M_3(|v_{i,j}|)\big)\!\mapsto_{\mathcal{A}}\!
\big( q_{i+1,j},M_2(|v_{i,j}|)\! +\! l_{i,j},M_3(|v_{i,j}|)\! +\! l'_{i,j}\big)$$
and $|w_{i,j}|\! =\! |v_{i,j}|.2^{l_{i,j}}.3^{l'_{i,j}}$. Moreover, the state $q_{n_j,j}$ is a final state, $M_2(|w_{n_j-1,j}|)\! =\! 0$, and 
$M_3(|w_{n_j-1,j}|)\! =\! 0$.\medskip 

 In particular, $|v_{0,0}|\! =\! 1\! =\! 2^0.3^0$. We will prove, by induction on $i\! <\! n_0$, that 
$$|w_{i,0}|\! =\! 2^{M_2(|w_{i,0}|)}.3^{M_3(|w_{i,0}|)}\mbox{,}$$ 
and $|w_{i,0}|\! =\! |v_{i+1,0}|$ if $i\! <\! n_0\! -\! 1$. Moreover, setting $c_0^i\! =\! M_2(|v_{i,0}|)$ and $c_1^i\! =\! M_3(|v_{i,0}|)$, we will prove that, for each $i\! <\! n_0\! -\! 1$, 
$a_{i,0}\! :\! (q_{i,0}, c_0^i,c_1^i)\!\mapsto_{\mathcal{A}}\! (q_{i+1,0},c_0^{i+1},c_1^{i+1})$.\medskip

 We have already seen that $|v_{0,0}|\! =\! 1\! =\! 2^0.3^0$. By hypothesis we can find a state $q_{1,0}\!\in\! Q$ and integers 
$l_{0,0},l'_{0,0}\!\in\!\{-1,0,1\}$ such that 
$$a_{0,0}\! :\!\big( q_0,M_2(|v_{0,0}|),M_3(|v_{0,0}|)\big)\!\mapsto_{\mathcal{A}}\! 
\big( q_{1,0},M_2(|v_{0,0}|)\! +\! l_{0,0},M_3(|v_{0,0}|)\! +\! l'_{0,0}\big)\mbox{,}$$
i.e., $a_{0,0}\! :\! (q_0,0,0)\!\mapsto_{\mathcal{A}}\! (q_{1,0},l_{0,0},l'_{0,0})$. Then 
$|w_{0,0}|\! =\! |v_{0,0}|.2^{l_{0,0}}.3^{l'_{0,0}}\! =\! 2^{l_{0,0}}.3^{l'_{0,0}}$. Now note that 
$|w_{0,0}.z_{0, 0}|\! =\! |u_{1,0}.v_{1,0}|\! =\! {\bf 0}^{N.6^l}$ and $|u_{1,0}|\! =\! |z_{0,0}|$. Thus 
$|v_{1,0}|\! =\! |w_{0,0}|\! =\! 2^{l_{0,0}}.3^{l'_{0,0}}$. Setting $c_0^0\! :=\! 0$, $c_1^0\! :=\! 0$, 
$c_0^1\! :=\! l_{0,0}\! :=\! M_2(|v_{1,0}|)$ and $c_1^1\! :=\! l'_{0,0}\! :=\! M_3(|v_{1,0}|)$, it holds that 
$$a_{0,0}\! :\! (q_0,c_0^0,c_1^0)\!\mapsto_{\mathcal{A}}\! (q_{1,0},c_0^1, c_1^1).$$
Assume now that, for each $i\! <\! n_0\! -\! 1$, it holds that 
$|w_{i,0}|\! =\! |v_{i+1,0}|\! =\! 2^{M_2(|w_{i,0}|)}.3^{M_3(|w_{i,0}|)}$ and 
$a_{i,0}\! :\! (q_{i,0},c_0^i,c_1^i)\!\mapsto_{\mathcal{A}}\! (q_{i+1,0}, c_0^{i+1},c_1^{i+1})$. We know that we can find a state $q_{n_0,0}\!\in\! Q$ and integers $l_{n_0-1,0},l'_{n_0-1,0}\!\in\!\{-1,0,1\}$ such that\medskip
 
\leftline{$a_{n_0-1,0}\! :\!\big( q_{n_0-1,0},M_2(|v_{n_0-1,0}|),M_3(|v_{n_0-1,0}|)\big)\!\mapsto_{\mathcal{A}}$}\smallskip
 
\rightline{$\big( q_{n_0,0},M_2(|v_{n_0-1,0}|)\! +\! l_{n_0-1,0},M_3(|v_{n_0-1,0}|)\! +\! l'_{n_0-1,0}\big) ,$}\medskip

\noindent i.e., $a_{n_0-1,0}\! :\! (q_{n_0-1,0},c_0^{n_0-1},c_1^{n_0-1})\!\mapsto_{\mathcal{A}}\! 
(q_{n_0},c_0^{n_0-1}\! +\! l_{n_0-1,0},c_1^{n_0-1}\! +\! l'_{n_0-1,0})$.\medskip

 Then $|w_{n_0-1, 0}|\! =\! |v_{n_0-1,0}|.2^{l_{n_0-1,0}}.3^{l'_{n_0-1,0}}\! =\! 
2^{c_0^{n_0-1}\! +\! l_{n_0-1,0}}.3^{c_1^{n_0-1}\! +\! l'_{n_0-1,0}}\! =\! 1$ since moreover, by hypothesis, 
$M_2(|w_{n_0-1,0}|)\! =\! M_3(|w_{n_0-1,0}|)\! =\! 0$.\medskip

 Finally we inductively proved the announced claim, and this shows that $a_{0,0}a_{1,0}\ldots a_{n_0-1,0}$ is accepted by 
$\mathcal{A}$, by final states and empty stack. On the other hand, $|w_{n_0-1,0}.z_{n_0-1,0}|\! =\! N.6^{l+n_0-1}$, 
$|w_{n_0-1,0}|\! =\! 1$ and $|u_{n_0}|\! =\! |z_{n_0-1,0}|$, thus ${|u_{n_0}|\! =\! N.6^{l+n_0-1}\! -\! 1}$. As above, we can prove, by induction on $j$, that, for every $j\!\in\!\om$, the finite word $a_{0,j}a_{1,j}\ldots a_{n_j-1,j}$ is in $L(\mathcal{A})$, and thus 
$\sigma$ is in $L(\mathcal{A})^\infty$.\medskip

 Conversely it is easy to see that if $\sigma\!\in\! L(\mathcal{A})^\infty$, then $g_{N,l}(\sigma )\!\in\!\mathcal{L}^\infty$.
\hfill{$\diamond$}\medskip




 We now come back to the proof of Proposition \ref{21}. By Claim 1, $\mathcal{L}$ is accepted by a one-counter automaton 
$\mathcal{C}$, and  there are at most 5 consecutive $\lambda$-transitions during a run of $\mathcal{C}$ on a finite word $w$.\medskip

 The alphabet $Y$ will be $\Si\cup\{ {\bf 0},{\bf 1},{\bf 2}\}$. We first define the language $B'$ which will be of the form $\mu\cup\mathcal{L}$, where 
$\mu$ is a finitary language over $Y$. We will moreover ensure that $\mu$ is accepted by a one-counter automaton, (for which   there are also at most 5 consecutive $\lambda$-transitions during a run  on a finite word) by finite states and empty stack, so that it will also be the case of $B'$, by non-determinism. The set $\mathcal{L}^\infty$ will look like $A^\infty$ on some compact set $K_{1,0}$. We actually define, for any natural numbers $N,l$ such that  $6$ does not  divide $N\!\geq\! 1$, some compact sets $K_{N,l}$. On the $K_{N,l}$'s, we will be able to control the complexity of $B'^\infty$, which will essentially be that of $\mathcal{L}^\infty$, and $\mathcal{L}^\infty$ will be complex. Out of the 
$K_{N,l}$'s, we do not know the complexity of $\mathcal{L}^\infty$. This is the reason why we introduce $\mu$. The set $\mu^\infty$ will be simple, will hide the possible complexity of $\mathcal{L}^\infty$ out of the union of the $K_{N,l}$'s, and will not hide the complexity of $\mathcal{L}^\infty$ on the $K_{N,l}$'s. We set $K_{N,l}\! :=\! g_{N,l}[\Si^\om]$. As $g_{N,l}$ is a homeomorphism onto its range, $K_{N,l}$ is compact. By Claim 2, $g_{N,l}^{-1}(\mathcal{L}^\infty\cap K_{N,l})\! =\! A^\infty$ for each $l\!\in\!\omega$.\medskip

 We are ready to define $\mu$. We set\medskip
 
\leftline{$\mu\! :=\!\{\ w\!\in\! Y^{<\omega}\mid\exists n\!\geq\! 2~~\exists (a_i)_{i<n}\!\in\!\Si^n~~
\exists (P_i)_{i<n}\!\in\! (\om\!\setminus\!\{ 0\} )^n~~\exists (Q_i)_{i<n}\!\in\!\om^n$}\smallskip

\rightline{$w\! =\! {\bf 0}~({{}^\frown}_{i<n}~a_i~{\bf 1}~{\bf 0}^{P_i}~{\bf 2}~{\bf 0}^{Q_i})~\wedge ~
(P_{n-2}\!\not=\! Q_{n-2}~\vee ~P_{n-1}\!\not=\! 6.P_{n-2})~\} .$}\medskip

 Note that all the words in $B'$ have the same form ${\bf 0}~({{}^\frown}_{i<n}~a_i~{\bf 1}~{\bf 0}^{P_i}~{\bf 2}~{\bf 0}^{Q_i})$. Note also that any finite concatenation of words of this form still has this form. We set 
$$S\! :=\!\{\ {\bf 0}~({{}^\frown}_{i\in\om}~a_i~{\bf 1}~{\bf 0}^{P_i}~{\bf 2}~{\bf 0}^{Q_i})\mid (a_i)_{i\in\om}\!\in\!\Si^\om 
~\wedge ~(P_i)_{i\in\om}\!\in\! (\om\!\setminus\!\{ 0\} )^\om ~\wedge ~(Q_i)_{i\in\om}\!\in\!\om^\om\ \}.$$
We now show that $\mu^\infty$ is ``simple". Note that 
$$\mu^\infty\! =\!\{\ \gamma\!\in\! Y^\omega\mid\forall l\!\in\!\omega\ \ \exists t\!\in\!\mu^l~\wedge ~{^\frown}_{i<l}\ t(i)\!\subseteq\!\gamma\ \} .$$ 
This shows that $\mu^\infty\!\in\!\bormtwo (Y^\omega )$.\medskip

 We first prove the result for $B'$ and  the class ${\bf \Si}_n^0$. Note that $B'^\infty\cap K_{1,0}\! =\!\mathcal{L}^\infty\cap K_{1,0}$ is not a 
${\bf\Pi}^0_n$ subset of $K_{1,0}$  since $g_{1,0}^{-1}(\mathcal{L}^\infty\cap K_{1,0})\! =\! A^\infty$, and $A^\infty$ is  ${\bf \Si}_n^0$-complete and hence not in the class  ${\bf\Pi}^0_n$,      so that $B'^\infty$ is not a 
${\bf\Pi}^0_n$ subset of $Y^\omega$. By 22.10 in \cite{Kechris94}, it remains to see that $B'^\infty$ is a ${\bf\Sigma}^0_n$-subset of 
$Y^\omega$.\medskip

 We define $F\! :\! S\!\setminus\!\mu^\infty\!\ra\! (\{\lambda\}\cup\mu )\!\times\!\mathcal{P}$ as follows. Let 
${\gamma\! :=\! {\bf 0}~({{}^\frown}_{i\in\om}~a_i~{\bf 1}~{\bf 0}^{P_i}~{\bf 2}~{\bf 0}^{Q_i})\!\in\! S\!\setminus\!\mu^\infty}$, and 
$(N,l)\!\in\!\mathcal{P}$ with $P_0\! =\! N.6^l$. If 
${\gamma\!\in\! K_{N,l}}$, then we put $F(\gamma )\! :=\! (\lambda ,N,l)$. If ${\gamma\!\notin\! K_{N,l}}$, then there is $i_0\!\in\!\om$ maximal for which $P_{i_0}\!\not=\! Q_{i_0}$ or $P_{i_0+1}\!\not=\! 6.P_{i_0}$. Let $(N',l')\!\in\!\mathcal{P}$ with $P_{i_0+1}\! =\! N'.6^{l'}$. We put $F(\gamma )\! :=\!\big( {\bf 0}~({{}^\frown}_{i\leq i_0}~a_i~{\bf 1}~{\bf 0}^{P_i}~
{\bf 2}~{\bf 0}^{Q_i})~a_{i_0+1}~{\bf 1}~{\bf 0}^{P_{i_0+1}}~{\bf 2}~{\bf 0}^{Q_{i_0+1}-1},N',l'\! +\! 1\big)$. We then set 
$R\! :=\! F[S\!\setminus\!\mu^\infty ]$.\medskip 

 Assume that $\gamma\!\in\! B'^\infty\!\setminus\!\mu^\infty$. Note that $\gamma\!\in\! S\!\setminus\!\mu^\infty$, so that 
$(t,N,l)\! :=\! F(\gamma )$ is defined, $t\!\subseteq\!\gamma$ and $\gamma\! -\! t\!\in\! K_{N,l}$. We define, for 
$(t,N,l)\!\in\! R$, $P_{t,N,l}\! :=\!\{\gamma\!\in\! Y^\om\mid t\!\subseteq\!\gamma ~\wedge ~\gamma\! -\! t\!\in\! K_{N,l}\}$ and 
$A_{t,N,l}\! :=\!\{\gamma\!\in\! P_{t,N,l}\mid\gamma\! -\! t\!\in\!\mathcal{L}^\infty\cap K_{N,l}\}$. Note that $P_{t,N,l}$ is compact, contained in $S\!\setminus\!\mu^\infty$, and $F(\gamma )\! =\! (t,N,l)$ if $\gamma\!\in\! P_{t,N,l}$. This shows that the $P_{t,N,l}$'s are pairwise disjoint and disjoint from $\mu^\infty$. Note also that $A_{t,N,l}$ is ${\bf \Si}_n^0$. The previous discussion shows that 
$B'^\infty\! =\!\mu^\infty\cup\bigcup_{(t,N,l)\in R}~A_{t,N,l}$, so that $B'^\infty$ is also in ${\bf \Si}_n^0$.\medskip

 For the class ${\bf \Pi}_n^0$, we note that 
$B'^\infty\! =\!\mu^\infty\!\setminus\! (\bigcup_{(t,N,l)\in R}~P_{t,N,l})\cup
\bigcup_{(t,N,l)\in R}~A_{t,N,l}\cap P_{t,N,l}$. Thus $\neg B'^\infty\! =\!\neg (\mu^\infty\cup\bigcup_{(t,N,l)\in R}~P_{t,N,l})\cup
\bigcup_{(t,N,l)\in R}~P_{t,N,l}\!\setminus\! A_{t,N,l}$. As $n\!\geq\! 3$, the first part is in ${\bf \Si}_n^0$, as well as the second, so that $B'^\infty$ is in ${\bf \Pi}_n^0$.\medskip
 
 To finish the proof, we first notice that it is easy to see that the finitary language $\mu$, as the finitary language $\mathcal{L}$,  is accepted by a non-deterministic one-counter automaton, for which   there are also at most 5 consecutive $\lambda$-transitions during a run  on a finite word,  by final states {\it and} empty counter.  Details are here left to the reader.  Then the language $B'\! =\!\mathcal{L}\cup\mu$ is also accepted by a  non-deterministic one-counter automaton,   by final states {\it and} empty counter, for which   there are also at most 5 consecutive $\lambda$-transitions during a run  on a finite word. In order to get the language $B$ from the language $B'$ we use a simple morphism which is a very particular case of a substitution.  If the alphabet of $B'$ is $Y\! :=\!\{ a_1,a_2,\ldots, a_k\}$ for some integer $k\!\geq\! 1$, then we add a letter $c$ to $Y$, set $Z\! :=\! Y\cup\{ c\}$, and consider the morphism $h\! :\! Y\!\ra\! Z^{<\om}$ defined by $h(a_i)\! =\! a_i c^6$ for each integer $i\!\in\! [1, k]$. This morphism is naturally extended to words and then to languages. Then we set $B\! =\! h(B')$. A word of $B$ is simply obtained from a word $w$ of $B'$ by adding $6$ letters $c$ after each letter of $w$. It is then easy to see that the language $B$ is accepted by a  non-deterministic {\it  real-time} one-counter automaton $\mathcal{A}$ by final states {\it and} empty counter. This automaton is simply obtained from a  non-deterministic one-counter automaton  $\mathcal{B}$, for which   there are also at most 5 consecutive $\lambda$-transitions during a run  on a finite word, accepting $B'$ by final states {\it and} empty counter. The simple idea is that the $\lambda$-transitions of $\mathcal{B}$ now occur during the reading by $\mathcal{A}$ of the letters $c$ in a word of $B$.\medskip
 
 Moreover it is easy to see that if $B'^\infty$ is ${\bf \Si}_n^0$-complete, (respectively,    ${\bf \Pi}_n^0$-complete), for some natural number $n\geq 3$, then $B^\infty$ is also 
${\bf \Si}_n^0$-complete, (respectively, ${\bf \Pi}_n^0$-complete).\medskip 
 
 Finally it is easy to use the morphism $f\! :\! Z\!\ra\!\{ {\bf 0},{\bf 1}\}^{<\om}$ defined by $f(a_j)\! =\! {\bf 0}^j{\bf 1}$ for every $j$ in 
$\{ 1,\ldots ,k\}$ and $f(c)\! =\! {\bf 0}^{k+1}{\bf 1}$. Then the language $f(B)\!\subseteq\!\{ {\bf 0},{\bf 1}\}^{<\om}$ is accepted by a non-deterministic {\it real-time} one-counter automaton by final states {\it and} empty counter, and it is easy to see that if $B^\infty$ is 
${\bf \Si}_n^0$-complete, (respectively, ${\bf \Pi}_n^0$-complete), for some natural number $n\!\geq\! 3$, then $f(B)^\infty$ is also 
${\bf \Si}_n^0$-complete, (respectively, ${\bf \Pi}_n^0$-complete).\ep\medskip
 
 We now finish the proof of the main result of this section.\medskip
  
\noindent\bf Proof of Theorem \ref{main}.(a).\rm\medskip

 Theorem \ref{pin} and the discussion after it provide, for $n\! =\! 1,2$, a regular finitary language $P_n$ such that $P_n^\infty$ is ${\bf\Pi}_n^0$-complete. So we are done if $n\!\leq\! 2$.\medskip 

This discussion also provides a finite alphabet $\Si_{P_3}$ and a finitary language $P_3\!\subseteq\!\Si_{P_3}^{<\om}$, accepted by a real-time one-counter automaton, by final states {\it and} empty stack, such that $P_3^\infty$ is ${\bf\Pi}_3^0$-complete. By an argument similar to the one used in the last paragraph of the proof of Proposition \ref{21}, it is possible to get $\Si_{P_3}=\{ {\bf 0},{\bf 1}\}$ with the same property.\medskip

By Proposition \ref{sub}, if $h\! :\!\Si_{P_3}\!\ra\! 2^{(\Si_{P_3}\cup\{\tla\} )^{<\om}}$ is the substitution defined by $a\!\mapsto\! L_3a$, then the language $h(P_3)$ is in $OCL(2)$, and $h(P_3)$ is accepted by a real-time two-iterated counter automaton accepting words by final states {\it and} empty stack. The proof of Theorem \ref{pin} shows that $h(P_3)^\infty$ is ${\bf\Pi}_4^0$-complete. Proposition \ref{21} provides a finite alphabet 
$\Si_{P_4}=\{ {\bf 0},{\bf 1}\}$ and a finitary language $P_4\!\subseteq\!\Si_{P_4}^{<\om}$, accepted by a real-time one-counter automaton, by final states {\it and} empty stack, such that $P_4^\infty$ is ${\bf\Pi}_4^0$-complete. It remains to repeat this argument with $n\!\geq\! 4$ instead of $3$.
\ep\medskip

 We obtained an inductive construction of languages $P_n$ accepted by one counter automata such that $P_n^\infty$ is ${\bf\Pi}_n^0$-complete. We can argue slightly differently, as follows. First we can show that Proposition  \ref{21} is valid if we replace in the hypothesis ``a real-time two-iterated counter automaton $\mathcal{A}$''  by ``a real-time k-iterated counter automaton $\mathcal{A}$, for some integer $k\geq 2$'';  the proof of this extension of Proposition \ref{21} is very similar to the proof of Proposition \ref{21}, the idea being that we have in this case to code the content of $k$ counters, but the ideas and the constructions of the proof are very similar, details are here left to the reader.  Then 
Theorem \ref{main}{\rm .(a)} now follows from Theorem \ref{pin} and from this extension of Proposition \ref{21}. Notice that we  only gave a detailed proof of Proposition  \ref{21} in the case of $k=2$ because it is easier to exposit and this case contains all the fundamental ideas of the proof of the extended case of an integer $k\geq 2$. 

\section{$\!\!\!\!\!\!$ ${\bf \Si}_n^0$-complete $\om$-powers}\indent
  
 We want to find an alphabet $\Gamma$ and a context free language $A\!\subseteq\!\Gamma^{<\omega}$ such that 
$A^\infty$ is ${\bf\Sigma}^0_n$-complete.\medskip

\noindent\bf Notation.\rm\ We will consider the bijection $\mathcal{P}\! :\!\omega\!\rightarrow\omega^2$ obtained by taking the diagonals with constant sum $(0,0)$, then $(1,0),(0,1)$, then $(0,2),(1,1),(2,0)$, then $(3,0),(2,1),(1,2),(0,3)$, and so on alternatively down and up in the second coordinate. Formally, we define ${\mathcal M}\! :\!\omega\!\rightarrow\!\omega$ by 
$${\mathcal M}(n)\! :=\!\mbox{max}\{ q\!\in\!\omega\mid\frac{q(q+1)}{2}\!\leq\! n\} .$$

 We will consider the bijection $<.,.>:\!\omega^2\!\rightarrow\omega$ defined by 
$$<N,p>:=\!\!\left\{\!\!\!\!\!\!\!\!
\begin{array}{ll}
& \frac{(N+p)(N+p+1)}{2}\! +\! N\mbox{ if }N\! +\! p\mbox{ is even,}\cr\cr
& \frac{(N+p)(N+p+1)}{2}\! +\! p\mbox{ if }N\! +\! p\mbox{ is odd.}
\end{array}
\right.$$
Its inverse bijection $\mathcal{P}\! :\!\omega\!\rightarrow\omega^2$ is given by 
$${\mathcal P}(q):=\!\!\left\{\!\!\!\!\!\!\!\!
\begin{array}{ll}
& \Big( q\! -\!\frac{{\mathcal M}(q)({\mathcal M}(q)+1)}{2},{\mathcal M}(q)\! -\! q\! +\!\frac{{\mathcal M}(q)({\mathcal M}(q)+1)}{2}\Big)\mbox{ if }{\mathcal M}(q)\mbox{ is even,}\cr\cr
& \Big( {\mathcal M}(q)\! -\! q\! +\!\frac{{\mathcal M}(q)({\mathcal M}(q)+1)}{2},q\! -\!\frac{{\mathcal M}(q)({\mathcal M}(q)+1)}{2}\Big)\mbox{ if }{\mathcal M}(q)\mbox{ is odd.}
\end{array}
\right.$$
If $\alpha\!\in\! 2^\omega$ and $M\!\in\!\omega$, then we define the $M$'th vertical $(\alpha )_M\!\in\! 2^\omega$ of $\alpha$ by setting 
$$(\alpha )_M(p)\! :=\!\alpha (<M,p>)$$ 
if $p\!\in\!\omega$. We also define the odd part $(\alpha )^1\!\in\! 2^\omega$ of $\alpha$ by setting 
$(\alpha )^1(q)\! :=\!\alpha (2q\! +\! 1)$ if $q\!\in\!\omega$.\medskip

\noindent\bf Example.\rm\ By 23.A in \cite{Kechris94}, the set 
$\mathcal{S}_3\! :=\!\{\ x\!\in\! 2^{\omega^2}\mid\exists m\!\in\!\omega ~~\exists^\infty n\!\in\!\omega ~~x(m,n)\! =\! {\bf 0}\ \}$ of double binary sequences having a vertical with infinitely many zeros is $\borathree$-complete. Note that the set 
$\mathcal{S}\! :=\!\{\ \alpha\!\in\! 2^\omega\mid\exists N\!\in\!\omega ~~\exists^\infty q\!\in\!\omega ~~\big( (\alpha )_{2N}\big)^1(q)\! =\! {\bf 1}\ \}$ is also $\borathree$-complete. Indeed, its definition shows that it is $\borathree$, and the map 
$c\! :\! 2^{\omega^2}\!\rightarrow\! 2^\omega$ defined by 
$$c(x)(q)\! :=\!\!\left\{\!\!\!\!\!\!\!\!
\begin{array}{ll}
& {\bf 0}\mbox{ if }{\mathcal P}(q)(0)\mbox{ is odd or }{\mathcal P}(q)(1)\mbox{ is even,}\cr\cr
& {\bf 1}\! -\! x\big(\frac{{\mathcal P}(q)(0)}{2},\frac{{\mathcal P}(q)(1)-1}{2}\big)\mbox{ if }{\mathcal P}(q)(0)
\mbox{ is even and }
{\mathcal P}(q)(1)\mbox{ is odd,}
\end{array}
\right.$$ 
is continuous and satisfies $\mathcal{S}_3\! =\! c^{-1}(\mathcal{S})$.\medskip

 Note that we also have $\mathcal{S}\! = \!\{\ \alpha\!\in\! 2^\omega\mid\exists N\!\in\!\omega  ~~\big( (\alpha )_{2N}\big)^1\!\in\! (1^{<\om}{\bf 1})^\infty\ \}$.  More generally we will consider in the sequel the $\om$-language $\mathcal{S}\! := \!\{\ \alpha\!\in\! 2^\omega\mid\exists N\!\in\!\omega  ~\big( (\alpha )_{2N}\big)^1\!\in\! L^\infty\ \}$, where $L$ is a finitary language over the alphabet $2\! :=\!\{ {\bf 0},{\bf 1}\}$, such that the $\om$-power $L^\infty$ is in the class ${\bf\Delta}^0_{\xi +1}\!\setminus\! {\bf\Sigma}^0_\xi$, where $\xi\!\geq\! 2$ is a countable ordinal.\medskip

 We will be able to take $\Gamma\! =\! 4\! :=\!\{ {\bf 0},{\bf 1},{\bf 2},{\bf 3}\}$. The language $A$ will be made of two pieces: we will have $A\! :=\!\mu\cup\pi$. Informally, the set $\pi^\infty$ will look like $\mathcal{S}$ on some nice compact set $K_0$. We actually define, for any natural number $l$, some compact set $K_l$. On the $K_l$'s, we will be able to control the complexity of $A^\infty$, which will essentially be that of $\pi^\infty$, and $\pi^\infty$ will be non-${\bf\Pi}^0_{\xi +1}$. Out of the $K_l$'s, we do not know the complexity of $\pi^\infty$. This is the reason why we introduce $\mu$. The set $\mu^\infty$ will be simple, will hide the possible complexity of $\pi^\infty$ out of the union of the $K_l$'s, and will not hide the complexity of 
$\pi^\infty$ on the $K_l$'s.\medskip

\noindent\bf Notation.\rm\ We will sometimes view 2 or 3 as alphabets, and sometimes view them as letters. To make this distinction clear, we will use the boldface notation {\bf 0}, {\bf 1}, {\bf 2}, {\bf 3} for the letters, and the lightface notation 2, 3 otherwise. So we have $2\! =\!\{ {\bf 0},{\bf 1}\}$, $3\! =\!\{ {\bf 0},{\bf 1},{\bf 2}\}$, and $4\! =\!\{ {\bf 0},{\bf 1},{\bf 2},{\bf 3}\}$. We set 
$$K_0\! :=\!\{\ ({^\frown}_{j\in\omega}\ {\bf 2}\ s_{2j}\ {\bf 3}\ s_{2j+1})\!\in\! 4^\omega\mid
\forall m\!\in\!\omega ~~s_m\!\in\! 2^m\ \} .$$ 

 The idea is to view an element $\alpha$ of the Cantor space $2^\omega$ as the concatenation of the diagonal finite binary sequences $s_m$, with $s_m\!\in\! 2^m$, using the bijection $\mathcal{P}$. In $K_0$, we introduce some separators of the $s_m$'s, ${\bf 2}$ and 
${\bf 3}$ alternatively, so that $\alpha$ is now seen as an element of $4^\omega$. Similarly, we set 
$K_{l+1}\! :=\!\{\ ({^\frown}_{j\in\omega}\ {\bf 3}\ s_{2j+1}\ {\bf 2}\ s_{2j+2})\!\in\! 4^\omega\mid
\forall m\!\geq\! 1~~s_m\!\in\! 2^{2l+2+m}\ \}$, for each $l\!\in\!\omega$, erasing the first $2l\! +\! 3$ diagonal finite binary sequences appearing in the elements of $K_0$. As the map $\varphi_l\! :\! K_l\!\rightarrow\! 2^\omega$, defined by 
$\varphi_l(\gamma )\! :=\! {^\frown}_{j\in\omega}\ s_{2j+\varepsilon}^{(-1)^{\varepsilon}}s_{2j+\varepsilon +1}^{(-1)^{\varepsilon +1}}$, where $\varepsilon\!\in\! 2$ is $0$ exactly when $l\! =\! 0$, is a homeomorphism, $K_l$ is compact.\medskip

 We define $f\! :\! 2\!\rightarrow\! 2^{4^{<\omega}}$ by 
$$f(a)\! :=\!\{\ a~t~{\bf 3}~u~v~{\bf 2}~w\!\in\! 4^{<\omega}\mid t,u,v,w\!\in\! 2^{<\omega}\wedge
\vert t\vert\! =\!\vert u\vert\mbox{ is even }\wedge\vert v\vert\! =\!\vert w\vert\!\geq\! 3\mbox{ is odd }\} .$$ 
The language $\pi$ will be of the form $\pi_0\cup\pi_1$, the latter language $\pi_1$ depending on some fixed language 
$L\!\subseteq\! 2^{<\omega}$. We first set 
$$\pi_0\! :=\!\{\ ({^\frown}_{j\leq N}\ {\bf 2}~s_{2j}\ {\bf 3}~s_{2j+1})~{\bf 2}~a\!\in\! 4^{<\omega}\mid 
N\!\in\!\omega ~\wedge ~(\forall q\!\leq\! 2N\! +\! 1~~s_q\!\in\! 2^{<\omega})~\wedge ~s_0\! =\!\lambda ~\wedge ~a\!\in\! 2\ \} .$$ 
Fix then $L\!\subseteq\! 2^{<\omega}$. We set $\pi_1\! :=\! f(L)$, extending $f$ as in Definition 3. We then set 
$\pi\! :=\!\pi_0\cup\pi_1$.\medskip

 In order to simplify further notation, we set, for $N\!\in\!\omega$ and $(k_m)_{m\in\omega}\!\in\!\omega^\omega$ fixed and 
$p,q\!\in\!\omega$, $M_q\! :=\! 2N\! +\! 2q\! +\! 2$ and $S^q_p\! :=\!\Sigma_{p\leq m\leq q}~(k_m\! +\! 1)$.\medskip

 The next lemma expresses the fact that $\pi^\infty$ looks like $\mathcal{S}$ on $K_0$.

\begin{lem} \label{similarly} $\varphi_0[\pi^\infty\cap K_0]\! =\!\{\ \alpha\!\in\! 2^\omega\mid\exists N\!\in\!\omega ~~
\big( (\alpha )_{2N}\big)^1\!\in\! L^\infty\ \}$.\end{lem}

\noindent\bf Proof.\rm\ Let $\gamma\!\in\!\pi^\infty\cap K_0$, and $\alpha\! :=\!\varphi_0(\gamma )$. We can write 
$\gamma\! =\! {^\frown}_{m\in\omega}\ w_m\! =\! {^\frown}_{j\in\omega}\ {\bf 2}\ s_{2j}\ {\bf 3}\ s_{2j+1}$, where 
$w_m\!\in\!\pi\!\setminus\{\lambda\}$ and $s_k\!\in\! 2^k$. As the first coordinate of $\gamma$ is $\bf 2$, $w_0$ is of the form 
$({^\frown}_{j\leq N}\ {\bf 2}\ s_{2j}\ {\bf 3}~\ s_{2j+1})~{\bf 2}~a$, and $a\! =\! s_{2N+2}(0)$, which exists since 
$\vert s_{2N+2}\vert\! =\! 2N\! +\! 2$.\medskip

 As the first coordinate of $\gamma$ not in $2$ after $w_0$ is $\bf 3$, $w_1$ is not in $\pi_0$. Thus $w_1\!\in\!\pi_1\! =\! f(L)$ is of the form ${^\frown}_{j\leq k_1}\ a^1_j~t^1_j~{\bf 3}\ u^1_jv^1_j\ {\bf 2}\ w^1_j$. Inductively on $j\!\leq\! k_1$,we see that 
$\vert t^1_j\vert\! =\! 2N\! =\!\vert u^1_j\vert$, $u^1_jv^1_j\! =\! s_{M_j+1}$, 
$\vert v^1_j\vert\! =\! M_j\! +\! 1\! -\! 2N\! =\!\vert w^1_j\vert$, and $w^1_j\!\subseteq\! s_{M_j+2}$. Indeed, for 
$j\! =\! 0$, this comes from the facts that $a~a^1_0~t^1_0\! =\! s_{M_0}$, 
$\vert t^1_0\vert\! =\!\vert u^1_0\vert$ and $\vert v^1_0\vert\! =\!\vert w^1_0\vert$. If $j\! <\! k_1$, then this comes from the facts that $w^1_j~a^1_{j+1}~t^1_{j+1}\! =\! s_{M_j+2}$, $\vert t^1_{j+1}\vert\! =\!\vert u^1_{j+1}\vert$ and 
$\vert v^1_{j+1}\vert\! =\!\vert w^1_{j+1}\vert$.\medskip

 As the first coordinate of $\gamma$ not in $2$ after $w_1$ is $\bf 3$, $w_2$ is not in $\pi_0$. Thus $w_2$ is of the form 
${^\frown}_{j\leq k_2}\ a^2_j~t^2_j\ {\bf 3}\ u^2_j\ v^2_j\ {\bf 2}\ w^2_j$. Inductively on $j\!\leq\! k_2$, we see that 
$\vert t^2_j\vert\! =\! 2N\! =\!\vert u^2_j\vert$, $u^2_jv^2_j\! =\! s_{M_{k_1+1+j}+1}$, 
$\vert v^2_j\vert\! =\! M_{k_1+1+j}\! +\! 1\! -\! 2N\! =\!\vert w^2_j\vert$, and $w^2_j\!\subseteq\! s_{M_{k_1+1+j}+2}$. Indeed, for 
$j\! =\! 0$, this comes from the facts that $w^1_{k_1}~a^2_0~t^2_0\! =\! s_{M_{k_1+1}}$, 
$\vert t^2_0\vert\! =\!\vert u^2_0\vert$ and $\vert v^2_0\vert\! =\!\vert w^2_0\vert$. If $j\! <\! k_2$, then this comes from the facts that $w^2_j~a^2_{j+1}~t^2_{j+1}\! =\! s_{M_{k_1+1+j}+2}$, $\vert t^2_{j+1}\vert\! =\!\vert u^2_{j+1}\vert$ and 
$\vert v^2_{j+1}\vert\! =\!\vert w^2_{j+1}\vert$.\medskip

 If we continue like this, we find $(k_m)_{m\in\omega}$ such that
$$w_{m+1}\! =\! {^\frown}_{j\leq k_{m+1}}\ a^{m+1}_j\ t^{m+1}_j\ {\bf 3}\ u^{m+1}_j\ v^{m+1}_j\ {\bf 2}\ w^{m+1}_j\mbox{,}$$ 
$\vert t^{m+1}_j\vert\! =\! 2N$, $\vert u^{m+1}_j\vert\! =\! 2N$, 
$\vert v^{m+1}_j\vert\! =\! 2(j\! +\! 1+\!S^m_1)\! +\! 1\! =\!\vert w^{m+1}_j\vert$, 
$u^{m+1}_j\ v^{m+1}_j\! =\! s_{\left( M_{j+1+S^m_1}\right) -1}$ for each $m$, 
$w^{m+1}_j~a^{m+1}_{j+1}~t^{m+1}_{j+1}\! =\! s_{M_{j+1+S^m_1}}$ for each 
$j\! <\! k_{m+1}$, and $w^{m+1}_{k_{m+1}}~a^{m+2}_0~t^{m+2}_0\! =\! s_{M_{S^{m+1}_1}}$.

\vfill\eject

 Moreover, for each $m$, ${^\frown}_{j\leq k_{m+1}}\ a^{m+1}_j\!\in\! L$. Note that 
$$a^{m+1}_j\! =\! s_{M_{j+S^m_1}}(M_{j+S^m_1}\! -\! 2N\! -\! 1)\! =\! 
s_{M_{j+S^m_1}}^{-1}(2N)$$
and $\alpha (<2N,2q\! +\!\eta >)\! =\! s_{2N+2q+\eta +1}^{(-1)^\eta}(2N)$ if 
$q\!\in\!\omega$ and $\eta\!\in\! 2$. Thus $\big( (\alpha )_{2N}\big)^1\!\in\! L^\infty$, as desired.\medskip

 Conversely, assume that $\big( (\alpha )_{2N}\big)^1\!\in\! L^\infty$ for some $N\!\in\!\omega$. We set 
$s_0\! :=\!\lambda$ and, for $j\! =\! 2q\! +\!\eta\!\in\!\omega$, 
$$s_{j+1}\! :=\! (<\alpha (\frac{j(j+1)}{2}),\cdots ,\alpha (\frac{(j+1)(j+2)}{2}\! -\! 1)>)^{(-1)^{\eta +1}}\mbox{,}$$ 
so that $\gamma\! :=\! {^\frown}_{j\in\omega}\ {\bf 2}\ s_{2j}\ {\bf 3}\ s_{2j+1}$ satisfies 
$\gamma\!\in\! K_0$ and $\varphi_0(\gamma )\! =\!\alpha$. We set 
$$w_0\! :=\! ({^\frown}_{j\leq N}\ {\bf 2}\ s_{2j}\ {\bf 3}~s_{2j+1})\ {\bf 2}\ s_{2N+2}(0)\mbox{,}$$ 
so that $w_0\!\in\!\pi_0$ and $w_0\!\subseteq\!\gamma$. Let $(v_m)_{m\in\omega}\!\in\! L^\omega$ with 
$\big( (\alpha )_{2N}\big)^1\! =\! {^\frown}_{m\in\omega}~v_m$.\medskip

 We set $w_{m+1}\! :=\! {^\frown}_{\Sigma_{i<m}\vert v_i\vert\leq q<\Sigma_{i\leq m}\vert v_i\vert}\ 
\Big(\big( s_{M_q}\! -\! s_{M_q}\vert (2q\! +\! 1)\big)\ {\bf 3}\ s_{M_q+1}\ {\bf 2}~
s_{M_q+2}\vert (2q\! +\! 3)\Big)$. Note that $w_{m+1}\!\in\! f(v_m)\!\subseteq\! f(L)\! =\!\pi_1$ since, with $j\! <\!\vert v_m\vert$ and 
$q\! :=\!\Sigma_{i<m}~\vert v_i\vert\! +\! j$, 
$$s_{M_q}(2q\! +\! 1)\! =\! s_{M_q}^{-1}(2N)\! =\!\alpha (<2N,2q\! +\! 1>)\! =\!\big( (\alpha )_{2N}\big)^1(q)
\! =\! v_m(j)$$ 
and $\gamma\! =\! {^\frown}_{m\in\omega}\ w_m\!\in\!\pi^\infty\cap K_0$, so that 
$\alpha\!\in\!\varphi_0[\pi^\infty\cap K_0]$ as desired.\hfill{$\square$}\medskip

\noindent\bf Notation.\rm\ We are ready to define $\mu$. The idea is that an infinite sequence 
containing a word in $\mu$ cannot be in $K_0$. We set $\mu\! :=\!\bigcup_{i\leq 2}~\mu_i$, where 
$\mu_0\! :=\!\{\ w\!\in\! 4^{<\omega}\mid\exists v\!\in\! 4^{<\omega}\!\setminus\!\{\lambda\}~~
v~{\bf 2}~{\bf 3}\!\subseteq\! w\ \}$, 
$$\mu_1\! :=\!\{\ w\!\in\! 4^{<\omega}\mid\exists v\!\in\! 4^{<\omega}~~
\exists v',v''\!\in\! 2^{<\omega}~~v~{\bf 3}~v'~{\bf 2}~v''~{\bf 3}\!\subseteq\! w
~\wedge ~\vert v''\vert\!\not=\!\vert v'\vert\! +\! 1\ \}\mbox{,}$$
and $\mu_2\! :=\!\{\ w\!\in\! 4^{<\omega}\mid\exists v\!\in\! 4^{<\omega}~~
\exists v',v''\!\in\! 2^{<\omega}~~v~{\bf 2}~v'~{\bf 3}~v''~{\bf 2}\!\subseteq\! w~\wedge ~
\vert v''\vert\!\not=\!\vert v'\vert\! +\! 1\ \}$. We now show that $\mu^\infty$ is ``simple". Note that 
$\mu^\infty\! =\!\{\ \gamma\!\in\! 4^\omega\mid\forall l\!\in\!\omega\ \ \exists t\!\in\!\mu^{<\omega}\ \ 
\vert t\vert\!\geq\! l~\wedge ~{^\frown}_{i<\vert t\vert}\ t(i)\!\subseteq\!\gamma\ \}$. This shows that 
$\mu^\infty\!\in\!\bormtwo (4^\omega )$.

\begin{thm}\label{complete} 
Let $\xi\!\geq\! 2$ be a countable ordinal. If 
$L^\infty\!\in\! {\bf\Delta}^0_{\xi +1}\!\setminus\! {\bf\Sigma}^0_\xi$, then $A^\infty$ is ${\bf\Sigma}^0_{\xi +1}$-complete.
\end{thm}

\noindent\bf Proof.\rm\ It is straightforward to prove that if $T\!\subseteq\! 2^\omega$ is ${\bf\Delta}^0_{\xi +1}\!\setminus\! {\bf\Sigma}^0_\xi$, then  
$\{\alpha\!\in\! 2^\om\mid\exists N\!\in\!\om ~~(\alpha )_N\!\in\! T\}$ is ${\bf\Sigma}^0_{\xi +1}$-complete. This, 22.10 in  \cite{Kechris94} and Lemma \ref{similarly} imply that $A^\infty\cap K_0\! =\!\pi^\infty\cap K_0$ is not a ${\bf\Pi}^0_{\xi +1}$ subset of $K_0$. Thus $A^\infty$ is not a ${\bf\Pi}^0_{\xi +1}$ subset of $4^\omega$. By 22.10 in \cite{Kechris94} again, it remains to see that $A^\infty$ is a ${\bf\Sigma}^0_{\xi +1}$ subset of $4^\omega$. We set, for $i,N\!\in\!\omega$ and $v\!\in\! 2^{2N+1}$, 
$$P_{v,i}\! :=\!\Big\{\ \alpha\!\in\! 2^\omega\mid {\bf 0}^{\frac{(M_i-1)M_i}{2}+2i+1}v\!\subseteq\!\alpha 
~\wedge ~\big( (\alpha )_{2N}\big)^1\! -\!\Big(\big( (\alpha )_{2N}\big)^1\vert i\Big)\!\in\! L^\infty\ \Big\}\mbox{,}$$
$$K_{v,i}\! :=\!\{\ \gamma\!\in\! 4^\omega\mid v\!\subseteq\!\gamma ~\wedge ~
\gamma\! -\! v\!\in\! K_{N+i+1}\ \} .$$
In other words, $P_{v,i}$ is the set of elements of the Cantor space starting with $M_i\! -\! 1$ diagonal finite binary sequences with only zeros, whose next diagonal starts with $2i\! +\! 1$ zeros, and such that the odd part of the $(2N)^{\mbox{th}}$ vertical, minus its initial segment of length $i$, is in $L^\infty$.\medskip

 The next claim, in the style of Lemma \ref{similarly}, essentially says that $\pi^\infty$ looks like $P_{v,i}$ on the compact set 
$K_{v,i}$.\medskip

\noi\bf Claim 1.\it\ Let $i,N\!\in\!\omega$ and $v\!\in\! 2^{2N+1}$. Then
$$\pi^\infty\cap K_{v,i}\! =\!\big\{\ \gamma\!\in\! 4^\omega\mid 
\delta\! :=\! ({^\frown}_{j\leq N+i}\ {\bf 2}\ {\bf 0}^{2j}\ {\bf 3}\ {\bf 0}^{2j+1})\ {\bf 2}\ {\bf 0}^{2i+1}\gamma\!\in\! K_0\wedge\varphi_0(\delta )\!\in\! P_{v,i}\ \big\} .$$\rm

 Indeed, let $\gamma\!\in\!\pi^\infty\cap K_{v,i}$, and 
$\delta\! :=\! ({^\frown}_{j\leq N+i}\ {\bf 2}\ {\bf 0}^{2j}\ {\bf 3}\ {\bf 0}^{2j+1})\ {\bf 2}\ {\bf 0}^{2i+1}\gamma$. Note that 
$\delta\!\in\! K_0$, so that $\alpha\! :=\!\varphi_0(\delta )$ is defined and starts with ${\bf 0}^{\frac{(M_i-1)M_i}{2}+2i+1}v$. We can write 
$$\gamma\! =\! v~{^\frown}_{j\in\omega}\ {\bf 3}\ s_{M_{i+j}+1}\ {\bf 2}\ s_{M_{i+j}+2}\! =\! 
{^\frown}_{m\in\omega}\ w_m\mbox{,}$$ where $s_k\!\in\! 2^k$ and 
$w_m\!\in\!\pi$. As the first coordinate of $\gamma$ not in $2$ is ${\bf 3}$, $w_0$ is of the form 
$${^\frown}_{j\leq k_0}\ a^0_j~t^0_j~{\bf 3}\ u^0_j\ v^0_j\ {\bf 2}\ w^0_j.$$ 
Inductively on $j\!\leq\! k_0$,we see that $\vert t^0_j\vert\! =\! 2N\! =\!\vert u^0_j\vert$, $u^0_jv^0_j\! =\! s_{M_{i+j}+1}$, 
$\vert v^0_j\vert\! =\! M_{i+j}\! +\! 1\! -\! 2N\! =\!\vert w^0_j\vert$, and $w^0_j\!\subseteq\! s_{M_{i+j}+2}$. Indeed, for 
$j\! =\! 0$, this comes from the facts that $a^0_0~t^0_0\! =\! v$, $\vert t^0_0\vert\! =\!\vert u^0_0\vert$ and 
$\vert v^0_0\vert\! =\!\vert w^0_0\vert$. If $j\! <\! k_0$, then this comes from the facts that 
$w^0_j~a^0_{j+1}~t^0_{j+1}\! =\! s_{M_{i+j}+2}$, $\vert t^0_{j+1}\vert\! =\!\vert u^0_{j+1}\vert$ and 
$\vert v^0_{j+1}\vert\! =\!\vert w^0_{j+1}\vert$.\medskip

 We then argue as in the proof of Lemma \ref{similarly} to get 
$(k_m)_{m\in\omega}$ with 
$$w_m\! =\! {^\frown}_{j\leq k_m}\ a^m_j\ t^m_j\ {\bf 3}\ u^m_j\ v^m_j\ {\bf 2}\ w^m_j\mbox{,}$$ 
$\vert t^m_j\vert\! =\! 2N$, $\vert u^m_j\vert\! =\! 2N$, 
$\vert v^m_j\vert\! =\! 2( i\! +\! j\! +\! 1+\!S^{m-1}_0)\! +\! 1\! =\!\vert w^m_j\vert$, 
$u^m_j\ v^m_j\! =\! s_{\left( M_{i+j+1+S^{m-1}_0}\right) -1}$ for each $m$, 
$w^m_j~a^m_{j+1}~t^m_{j+1}\! =\! s_{M_{i+j+1+S^{m-1}_0}}$ for each $j\! <\! k_m$, and 
$w^m_{k_m}~a^{m+1}_0~t^{m+1}_0\! =\! s_{M_{i+S^m_0}}$. Moreover, for each $m$, 
${^\frown}_{j\leq k_m}\ a^m_j\!\in\! L$. Note that 
$a^m_j\! =\! s_{M_{i+j+S^{m-1}_0}}(M_{i+j+S^{m-1}_0}\! -\! 2N\! -\! 1)\! =\! 
s_{M_{i+j+S^{m-1}_0}}^{-1}(2N)$. Thus 
$$\begin{array}{ll}
\big( (\alpha )_{2N}\big)^1\! -\!\Big(\big( (\alpha )_{2N}\big)^1\vert i\Big)\!\!\!\!
& =\!\Big(\big( (\alpha )_{2N}\big)^1(i),\big( (\alpha )_{2N}\big)^1(i\! +\! 1),\cdots\Big)\cr
& =\!\big(\alpha (<2N,2i\! +\! 1>),\alpha (<2N,2i\! +\! 3>),\cdots\big)\cr
& =\! (s_{M_i}^{-1}(2N),s_{M_i+2}^{-1}(2N),\cdots )\cr
& =\! (a^0_0,a^0_1,\cdots ,a^0_{k_0},a^1_0,\cdots ,a^1_{k_1},\cdots )
\end{array}$$
is in $L^\infty$, as desired.\medskip

 Conversely, assume that $\gamma\!\in\! 4^\omega$, 
$\delta\! :=\! ({^\frown}_{j\leq N+i}\ {\bf 2}\ {\bf 0}^{2j}\ {\bf 3}\ {\bf 0}^{2j+1})\ {\bf 2}\ {\bf 0}^{2i+1}\gamma\!\in\! K_0$, and 
$\alpha\! :=\!\varphi_0(\delta )$ is in $P_{v,i}$. Then $\gamma\!\in\! K_{v,i}$. We set, for $j\! =\! 2q\! +\!\eta\!\geq\! M_i$, 
$$s_{j+1}\! :=\! (<\alpha (\frac{j(j+1)}{2}),\cdots ,\alpha (\frac{(j+1)(j+2)}{2}\! -\! 1)>)^{(-1)^{\eta +1}}\mbox{,}$$ 
so that $\gamma\! =\! v~{^\frown}_{j\in\omega}\ {\bf 3}\ s_{M_{i+j}+1}\ {\bf 2}\ s_{M_{i+j}+2}$. Let 
$(v_m)_{m\in\omega}\!\in\! L^\omega$ with 
$$\big( (\alpha )_{2N}\big)^1\! -\!\Big(\big( (\alpha )_{2N}\big)^1\vert i\Big)\! =\! {^\frown}_{m\in\omega}~v_m.$$ 

\vfill\eject

 We set\medskip

\leftline{$w_0\! :=\! v~\big( {\bf 3}\ s_{M_i+1}\ {\bf 2}\ s_{M_i+2}\vert (2i\! +\! 3)\big)$}\smallskip

\rightline{$\bigg( {^\frown}_{i<q<i+\vert v_0\vert}\ 
\Big(\big( s_{M_q}\! -\! s_{M_q}\vert (2q\! +\! 1)\big)\ {\bf 3}\ s_{M_q+1}\ {\bf 2}~s_{M_q+2}\vert (2q\! +\! 3)\Big)\bigg)\mbox{,}$}\medskip
 
\noindent so that $w_0\!\in\!\pi_1$ and $w_0\!\subseteq\!\gamma$. We then set
$$w_{m+1}\! :=\! {^\frown}_{i+\Sigma_{k\leq m}\vert v_k\vert\leq q<i+\Sigma_{k\leq m+1}\vert v_k\vert}~
\Big(\big( s_{M_q}\! -\! s_{M_q}\vert (2q\! +\! 1)\big)\ {\bf 3}\ s_{M_q+1}\ {\bf 2}~s_{M_q+2}\vert (2q\! +\! 3)\Big) .$$
Note that $w_{m+1}\!\in\! f(v_{m+1})\!\subseteq\! f(L)\! =\!\pi_1$ since, with $j\! <\!\vert v_{m+1}\vert$ and 
$q\! :=\! i\! +\!\Sigma_{k\leq m}~\vert v_k\vert\! +\! j$, 
$$s_{M_q}(2q\! +\! 1)\! =\! s_{M_q}^{-1}(2N)\! =\!\alpha (<2N,2q\! +\! 1>)\! =\!\big( (\alpha )_{2N}\big)^1(q)
\! =\! v_{m+1}(j)$$ 
and $\gamma\! =\! {^\frown}_{m\in\omega}\ w_m\!\in\!\pi^\infty$, as desired.\hfill{$\diamond$}\medskip

 The next claim provides a characterization of $A^\infty$ giving an upper bound on its topological complexity.\medskip

\noindent\bf Claim 2.\it\ Let $\gamma\!\in\! 4^\omega$. Then 
$$\gamma\!\in\! A^\infty ~\Leftrightarrow ~\gamma\!\in\!\mu^\infty ~\vee ~\gamma\!\in\!\pi^\infty\cap K_0\vee ~
\exists t\!\in\!\{\lambda\}\cup\mu~~(t\!\subseteq\!\gamma ~\wedge ~\exists i,N\!\in\!\omega ~~
\exists v\!\in\! 2^{2N+1}~~\gamma\! -\! t\!\in\!\pi^\infty\cap K_{v,i}).$$\rm

 Indeed, the right to left implication is clear. So assume that $\gamma\!\in\! A^\infty\!\setminus\!\mu^\infty$. Note that we can find $(v_j)_{j\in\omega}\!\in\! (2^{<\omega})^\omega$, $(a_j)_{j\in\omega}\!\in\!\{ {\bf 2},{\bf 3}\}^\omega$ and 
$(w_m)_{m\in\omega}\!\in\! A^\omega$ with 
$\gamma\! =\! {^\frown}_{j\in\omega}\ v_j\ a_j\! =\! {^\frown}_{m\in\omega}\ w_m$. As $\gamma\!\notin\!\mu^\infty$, there is $m_0\!\in\!\omega$ such that $w_m$ is of the form 
${^\frown}_{j\leq k_m}\ a^m_j\ t^m_j~{\bf 3}\ u^m_j\ v^m_j\ {\bf 2}\ w^m_j$ if $m\!\geq\! m_0$. Moreover, we may assume that $\vert t^m_j\vert\! =\!\vert t^{m_0}_0\vert$ is even and 
$\vert w^m_j\vert\! =\!\vert w^{m_0}_0\vert\! +\! 2\big(S^{m-1}_{m_0}\! +\! j\big)\!\geq\! 3$ is odd if 
$m\!\geq\! m_0$, and that $m_0$ is minimal with these properties.\medskip

\noindent\bf Case 1.\rm\ $m_0\! =\! 0$.\medskip

 We set $t\! :=\!\lambda$, $i\! :=\!\frac{\vert w^0_0\vert -3}{2}$, $N\! :=\!\frac{\vert t^0_0\vert}{2}$, 
$v\! :=\! a^0_0~t^0_0$ and $\delta\! :=\! {^\frown}_{m\geq m_0}\ w_m$, so that $\delta\!\in\!\pi^\infty\cap K_{v,i}$ and $\gamma\! =\!\delta$ is as desired.\medskip

\noindent\bf Case 2.\rm\ $\exists m\! <\! m_0$ such that $w_m\!\in\!\mu$.\medskip

 We set $t\! :=\! {^\frown}_{m<m_0}\ w_m$, $i\! :=\!\frac{\vert w^{m_0}_0\vert -3}{2}$, $N\! :=\!\frac{\vert t^{m_0}_0\vert}{2}$, 
$v\! :=\!a^{m_0}_0~t^{m_0}_0$ and $\delta\! :=\! {^\frown}_{m\geq m_0}\ w_m$, so that $t\!\in\!\mu$, $t\!\subseteq\!\gamma$, 
$\delta\!\in\!\pi^\infty\cap K_{v,i}$ and $\gamma\! -\! t\! =\!\delta$ is as desired.\medskip

\noindent\bf Case 3.\rm\ $\exists m\! <\! m_0$ such that $w_m$ is of the form 
$({^\frown}_{j\leq N_m}\ {\bf 2}~s_{2j}\ {\bf 3}~s_{2j+1})~{\bf 2}~a$.\medskip

 If $m\!\geq\! 1$, then $t\! :=\! {^\frown}_{m<m_0}\ w_m$ is in $\mu$, and we argue as in Case 2. So we may assume that $m\! =\! 0$. If $\gamma\!\in\! K_0$, then $\gamma\!\in\!\pi^\infty$. So we may assume that 
$\gamma\!\notin\! K_0$, which gives $j_0\!\in\!\omega$ such that 
$\vert v_{j_0+1}\vert\!\not=\!\vert v_{j_0}\vert\! +\! 1$. If $J\! >\! j_0$, then ${^\frown}_{j\leq J}\ v_j\ a_j\!\in\!\mu$. We choose $J$ big enough to ensure that ${^\frown}_{m<m_0}\ w_m\!\subseteq\! {^\frown}_{j\leq J}\ v_j\ a_j$. We then choose $m_1\!\geq\! m_0$ such that ${^\frown}_{j\leq J}\ v_j\ a_j\!\subseteq\! {^\frown}_{m<m_1}\ w_m$. We set 
$t\! :=\! {^\frown}_{m<m_1}\ w_m$, $i\! :=\!\frac{\vert w^{m_1}_0\vert -3}{2}$, $N\! :=\!\frac{\vert t^{m_1}_0\vert}{2}$, 
$v\! :=\! a^{m_1}_0~t^{m_1}_0$ and $\delta\! :=\! {^\frown}_{m\geq m_1}\ w_m$, so that $t\!\in\!\mu$, 
$t\!\subseteq\!\gamma$, $\delta\!\in\!\pi^\infty\cap K_{v,i}$ and $\gamma\! -\! t\! =\!\delta$ is as desired.

\vfill\eject

\noindent\bf Case 4.\rm\ $m_0\!\geq\! 1$ and $w_m$ is of the form 
${^\frown}_{j\leq k_m}\ a^m_j\ t^m_j~{\bf 3}\ u^m_j\ v^m_j\ {\bf 2}\ w^m_j$ if 
$m\! <\! m_0$.\medskip

 The minimality of $m_0$ gives $j\!\leq\! k_{m_0}$ such that $\vert t^{m_0-1}_0\vert\!\not=\!\vert t^{m_0}_j\vert$ or 
$$\vert w^{m_0}_j\vert\!\not=\!\vert w^{m_0-1}_0\vert\! +\! 2(k_{m_0-1}\! +\! 1\! +\! j).$$ 

 We set $t\! :=\! {^\frown}_{m\leq m_0}\ w_m$, $i\! :=\!\frac{\vert w^{m_0+1}_0\vert -3}{2}$, 
$N\! :=\!\frac{\vert t^{m_0+1}_0\vert}{2}$, $v\! :=\! a^{m_0+1}_0~t^{m_0+1}_0$ and 
$\delta\! :=\! {^\frown}_{m>m_0}\ w_m$, so that $t\!\in\!\mu$, $t\!\subseteq\!\gamma$, 
$\delta\!\in\!\pi^\infty\cap K_{v,i}$ and $\gamma\! -\! t\! =\!\delta$ is as desired.\hfill{$\diamond$}\medskip

 Note that $P_{v,i}$ is a ${\bf\Delta}^0_{\xi +1}$ subset of $2^\omega$. By Claim 1, $\pi^\infty\cap K_{v,i}$ is a 
${\bf\Delta}^0_{\xi +1}$ subset of $4^\omega$. By Claim 2, $A^\infty$ is a ${\bf\Sigma}^0_{\xi +1}$ subset of $4^\omega$.\hfill{$\square$}\medskip 
 
 It remains to see that $A$ is accepted by a one-counter automaton. We first check that $\mu_0,\mu_1,\mu_2,\pi_0$, $\pi_1$ are accepted by a one-counter automaton. The language $\mu_0$ is not only accepted by a one-counter automaton, it is in fact regular. 
 
\begin{lem} \label{mu0} The language $\mu_0\! :=\!\{\ w\!\in\! 4^{<\omega}\mid\exists v\!\in\! 4^{<\omega}\!\setminus\!\{\lambda\}~~
v~{\bf 2}~{\bf 3}\!\subseteq\! w\ \}$ is regular.\end{lem}

\proo  It is easy to construct  a finite  automaton accepting $\mu_0$. \hfill{$\square$}

\begin{lem} \label{pi0} The language 
$$\pi_0\! :=\!\{\ ({^\frown}_{j\leq N}\ {\bf 2}~s_{2j}\ {\bf 3}~s_{2j+1})~{\bf 2}~a\!\in\! 4^{<\omega}\mid 
N\!\in\!\omega ~\wedge ~(\forall q\!\leq\! 2N\! +\! 1~~s_q\!\in\! 2^{<\omega})~\wedge ~s_0\! =\!\lambda ~\wedge ~a\!\in\! 2\ \}$$ 
is regular.\end{lem}

\proo    It is again easy to construct  a finite  automaton accepting  the language  $\pi_0\!$. The details are here left to the reader.\hfill{$\square$}

\begin{lem}\label{mu12} 
The languages 
$$\begin{array}{ll}
& \mu_1\! :=\!\{\ w\!\in\! 4^{<\omega}\mid\exists v\!\in\! 4^{<\omega}~~\exists v',v''\!\in\! 2^{<\omega}~~v~{\bf 3}~v'~{\bf 2}~v''~{\bf 3}\!\subseteq\! w~\wedge ~\vert v''\vert\!\not=\!\vert v'\vert\! +\! 1\ \}\mbox{,}\cr\cr
& \mu_2\! :=\!\{\ w\!\in\! 4^{<\omega}\mid\exists v\!\in\! 4^{<\omega}~~\exists v',v''\!\in\! 2^{<\omega}~~v~{\bf 2}~v'~{\bf 3}~v''~{\bf 2}\!\subseteq\! w~\wedge ~\vert v''\vert\!\not=\!\vert v'\vert\! +\! 1\ \}
\end{array}$$ 
are accepted by real-time one-counter automata accepting words by final states {\it and} empty stack.
\end{lem}

\proo We indicate informally the idea of the construction of a real-time one-counter automaton $\mathcal{A}$ accepting the language $\mu_1\!$ by final states {\it and} empty stack. The automaton can use its finite control to check that the input word has an initial segment of the form ${\bf 3}~v'~{\bf 2}~v''~{\bf 3}\!$ for some finite words $v',v''\!\in\! 2^{<\omega}$. Moreover the automaton $\mathcal{A}$  can use its counter and the non-determinism to check that $\vert v''\vert\!\not=\!\vert v'\vert\! +\! 1$.\medskip

 If the automaton guesses that $\vert v''\vert\!>\!\vert v'\vert\! +\! 1$, then it increases its counter by $1$ for each letter of $v'$ and for the next letter ${\bf 2}$ which is read; next, while reading the segment $v''$, it decreases its counter by $1$ for each letter of $v''$ which is read, checking that the counter value becomes zero before ending the reading of  $v''$. On the other hand, if the automaton guesses that $\vert v''\vert\!<\!\vert v'\vert\! +\! 1$, then the automaton $\mathcal{A}$ begins to read a non null number $k$ of letters of $v'$ without increasing the counter,  guessing that $\vert v''\vert\!=\vert v'\vert\! +\! 1\! -\! k$; then it decreases the counter by $1$ for  each letter of $v'$ and for the next letter ${\bf 2}$ which is read; and the automaton checks that 
$\vert v''\vert\!=\vert v'\vert\! +\! 1\! -\! k$ by decreasing the counter by $1$ for each letter of $v''$ which is read.

\vfill\eject

 Similar ideas are used in the case of the language  $\mu_2$. The details are here left to the reader.\hfill{$\square$}
 
\begin{lem}\label{pi1} 
Let $L$ be a finitary language over $2$ accepted by a real-time one-counter automaton, by final states 
{\it and} empty stack. Then the language\medskip

\leftline{$\pi_1\! :=\!\{\ ({^\frown}_{j\leq k}\ t_j~{\bf 3}\ u_j\ v_j\ {\bf 2}\ w_j)\!\in\! 4^{<\omega}\mid k\!\in\!\omega ~
\wedge ~\vert t_j\vert\! =\!\vert u_j\vert\! +\! 1\mbox{ is odd }\wedge ~\vert v_j\vert\! =\!\vert w_j\vert\geq\! 3
\mbox{ is odd }~\wedge$}\smallskip

\rightline{${{}^\frown}_{j\leq k}~t_j(0)\!\in\! L\ \}$}\medskip

\noindent is in $OCL(2)$, and $\pi_1$ is accepted by a real-time two-iterated counter automaton, by final states {\it and} empty stack. If moreover $L$ is rational, hence accepted by a real-time finite automaton (without any counter)  by final states, then $\pi_1$  is in  $OCL(1)$ and is  accepted by a real-time one counter automaton, by final states {\it and} empty stack.
\end{lem}

\proo
Let $L$ be a finitary language over $2$ accepted by a real-time one-counter automaton $\mathcal{A}$, by final states {\it and} empty stack. We assume that the stack alphabet of 
$\mathcal{A}$ is equal to $\Gamma\! :=\!\{ Z_0,z_0\}$, and we informally explain the behaviour of a real-time two-iterated counter automaton $\mathcal{B}$ which will accept the language $\pi_1$ by final states {\it and} empty stack. The stack alphabet of $\mathcal{B}$ is equal to $\Gamma'\! :=\!\{ Z_0,z_0,z_1\}$ and the content of its  stack is always of the form $(z_1)^{n_1} (z_0)^{n_0}Z_0$ for some natural numbers $n_0,n_1$. The automaton  $\mathcal{B}$ can use its finite control to check that the input word is of the form 
$({^\frown}_{j\leq k}\ t_j~{\bf 3}\ u_j\ v_j\ {\bf 2}\ w_j)\!\in\! 4^{<\omega}$, for some $t_j, u_j, v_j, w_j\!\in\! 2^{<\omega}$.\medskip

 We now explain the behaviour of the automaton  $\mathcal{B}$ using its stack when reading a word of the form ${^\frown}_{j\leq k}\ t_j~{\bf 3}\ u_j\ v_j\ {\bf 2}\ w_j$. At the beginning the automaton reads $t_0(0)$ and it simulates the automaton  $\mathcal{A}$ with stack alphabet $\Gamma$. Then when reading the remaining part of $t_0$  it uses the stack letter $z_1$ and pushes a letter $z_1$ for each letter of $t_0$  read. When reading $u_0$ the automaton $\mathcal{B}$ pops a letter $z_1$ for each letter read until all letters $z_1$ have been popped from the stack. Again when reading $v_0$ the automaton pushes a letter $z_1$ for each letter of $v_0$  read and it pops a letter $z_1$ for each letter of $w_0$ read until all letters $z_1$ have been popped from the stack. The next letter to  be read is $t_1(0)$ and the automaton $\mathcal{B}$ simulates again the automaton $\mathcal{A}$ while reading this letter. Moreover it uses the ``second counter'' at the top of its stack with letters $z_1$ to check that 
$\vert t_1\vert\! =\!\vert u_1\vert\! +\! 1\wedge ~\vert v_1\vert\! =\!\vert w_1\vert$. The reading continues like that and the finite control can be used to check that $\vert t_j\vert\!$ is odd and $\!\vert w_j\vert\geq\! 3 $  is odd for every $j$. Finally after having read the letter $t_k(0)$ the automaton $\mathcal{B}$ has simulated the automaton $\mathcal{A}$ 
 on ${{}^\frown}_{j\leq k}~t_j(0)\!$ and it can check by final states and empty stack that ${{}^\frown}_{j\leq k}~t_j(0)\!\in\! L$; the automaton has only now to check the form of 
$t_k~{\bf 3}\ u_k\ v_k\ {\bf 2}\ w_k$, ending the reading in an accepting state and with an empty stack. This finishes the proof in the case of a language $L$ accepted by a real-time one-counter automaton $\mathcal{A}$, by final states {\it and} empty stack.\medskip

The proof is very similar and just simpler in the case of a language $L$ which is rational and accepted by a real-time finite automaton by final states. 
\hfill{$\square$}\medskip

\noindent {\bf Proof of Theorem 1.(b).}  The proof of Theorem 1.2 in \cite{Fink-Lec2} shows that if 
$$S_1\! =\!\{ w\!\in\! 2^{<\om}\mid {\bf 0}\!\subseteq\! w\vee\exists k\!\in\!\omega ~~{\bf 1}{\bf 0}^k{\bf 1}\!\subseteq\! w\}\mbox{,}$$ 
then $S_1^\infty$ is $\boraone$-complete. Note that $S_1$ is regular, and thus accepted by a one-counter automaton.\medskip

 Theorem 2 in  \cite{Fink-Lec2}  provides a finitary language $S_2$ which is accepted by a one-counter automaton and such that $S_2^\infty$ is ${\bf\Sigma}_2^0$-complete. So we are done if $n\!\leq\! 2$.\medskip
 
  Note that the language $L\! :=\!\{ w\!\in\! 2^{<\omega}\mid\exists j\! <\!\vert w\vert ~~w(j)\! =\! {\bf 1}\}$, the set of finite binary sequences having at least one coordinate equal to ${\bf 1}$, is rational and hence is  accepted by a real-time finite automaton, by final states. By Lemma \ref{pi1}, the language $\pi_1$ associated with $L$ is in $OCL(1)$, and $\pi_1$ is accepted by a real-time one counter automaton, by final states {\it and} empty stack.  By Lemmas \ref{mu0}, \ref{mu12}, \ref{pi0} and the non-determinism, this is also the case of $\mu\cup\pi$. By Theorem \ref{complete}, $(\mu\cup\pi )^\infty$ is $\borathree$-complete since 
$L^\infty\! =\!\mathbb{P}_\infty\!\in\!\bormtwo\!\setminus\!\boratwo\!\subseteq\!\borthree\!\setminus\!\boratwo$.\medskip
  
    The proof of Theorem 1 (a) provides a finite alphabet $\Si$ and a finitary language $P_3\!\subseteq\!\Si^{<\omega}$, accepted by a real-time one counter automaton, by final states {\it and} empty stack, such that $P_3^\infty$ is ${\bf\Pi}_3^0$-complete. Coding letters of $\Si$ with finite words over $2$ if necessary, we may assume that $\Si\! =\! 2$. By Lemma \ref{pi1}, the language $\pi_1$ associated with $P_3$ is in $OCL(2)$, and $\pi_1$ is accepted by a real-time two-iterated counter automaton, by final states {\it and} empty stack. By Lemmas \ref{mu0}, \ref{mu12}, \ref{pi0} and the non-determinism, this is also the case of $\mu\cup\pi$.  By Theorem \ref{complete}, $(\mu\cup\pi )^\infty$ is ${\bf \Sigma}_4^0$-complete since 
$(P_3)^\infty\! \!\in\! {\bf \Pi}_3^0\!\setminus\! {\bf \Sigma}_3^0\!\subseteq\! {\bf \Delta}_4^0\!\setminus\!{\bf \Sigma}_3^0$. Proposition \ref{21} provides a finite alphabet $\Si_{S_4}$ and a finitary language $S_4\!\subseteq\!\Si_{S_4}^{<\om}$, accepted by a real-time one-counter automaton, by final states {\it and} empty stack, such that $S_4^\infty$ is ${\bf\Sigma}_4^0$-complete.\medskip 
    
 It remains to repeat this argument with $n\!\geq\! 4$ instead of $3$.\hfill{$\square$}

\end{document}